\documentclass[11pt]{article}

\usepackage{amsmath,amsthm,verbatim,amssymb,amsfonts,amscd, graphicx}
\usepackage{graphics}
\theoremstyle{plain}

\newtheorem{thm}{Theorem}[section]
\newtheorem{question}[thm]{Question}

\newtheorem{lem}[thm]{Lemma}
\newtheorem{rem}[thm]{Remark}

\newtheorem{claim}[thm]{Claim}

\newtheorem{defn}[thm]{Definition}

\usepackage{extarrows}
\usepackage{mathrsfs}
\usepackage{galois}
\usepackage[colorlinks,citecolor = red, linkcolor=red,hyperindex,CJKbookmarks]{hyperref}
\usepackage[all]{hypcap}
\usepackage[all]{xy}
\usepackage{tikz-cd}

\begin{document}

\title{Accumulation point theorem for generalized log canonical thresholds}
\date{}
\setcounter{footnote}{-1}
\author{JIHAO LIU \footnote{\textit{Date}: 	Oct 30th, 2018.}}
\maketitle

\textbf{ABSTRACT}. In this paper we show that the set of accumulation points of generalized log canonical thresholds for certain DCC sets comes from the set of generalized log canonical thresholds of dimension $1$ less of the same DCC sets. 

\tableofcontents
\section{Introduction}

Log canonical thresholds are important invariants of singularities that play a fundamental role in higher dimensional birational geometry, e.g. \hyperref[Bir07]{[Bir07]}, \hyperref[HMX14]{[HMX14]}. Recently there has been substantial progress in understanding their behaviors with the proof of the ACC conjecture for log canonical thresholds:

\begin{thm}\label{0101}\rm (ACC for log canonical thresholds, \hyperref[HMX14]{[HMX14,Theorem 1.1]})
Let $d>0$ be an integer, $I\subset [0,1]$, $J\subset [0,\infty)$ be two DCC sets. Then 
\begin{center}
$LCT_d(I;J):=\{lct(X,B;D)|(X,B)$ is log canonical of dimension $d, B\in I, D\in J\}$ 
\end{center}
satisfies the ACC.

\end{thm}

In recent years, Birkar and Zhang has introduced developed the theory of \textit{generalized pairs}, which has  proved to be very useful in many cases. In \hyperref[BZ16]{[BZ16]},  the notion of \textit{generalized log canonical thresholds} is defined, and it is shown that they satisfy the ACC property:
\begin{thm}\label{0102}\rm (ACC of generalized log canonical thresholds, \hyperref[BZ16]{[BZ16, Theorem 1.5]})
Let $d>0$ be an integer, $I,J\subset [0,\infty)$ be two DCC sets. Then the set
\begin{center}
$GLCT_{d}(I;J):=\{glct(X',B'+M';C'+N')|(X',B'+M')$ is generalized LC of dimension $d$,

$B'\in I, C'\in J, M=\sum m_jM_j, N=\sum n_jN_j$ where each $M_j,N_j$ is nef Cartier, $ m_j\in I, n_j\in J\}$
\end{center}
satisfies the ACC.
\end{thm}

Notice that when $M$ and $N$ both equal to $0$, the generalized log canonical thresholds are exactly the usual log canonical thresholds; and when $J=\mathbb N^+$, $N'\not=0$, in some cases the generalized log canonical thresholds behave quite similarly to minimal log discrepancies. Thus, a deep understanding of generalized log canonical thresholds may help us greatly in the proof of ACC for minimal log discrepancies.

In this paper we focus on the accumulation point conjectures (theorems) associated to the three ACC conjectures above. That is, how much do we know about the accumulation points of $LCT_{d}(I,J)$, $GLCT_{d}(I,J)$. In particular, there are three natural questions:
\begin{question}\label{0103}\rm

Let $K=K(I,J,d)$ be one of the two sets above. 
 
(1) Under what conditions do we know that the accumulation points of $K$ are rational numbers?

(2) Under what conditions do we know that we can only take accumulation points finitely many times, and can we find a bound for this number? That is, we define $K:=K_0\subset K_1\subset\dots$ such that $K_n:=\bar K_{n-1}$ for every $n\geq 1$. Does there exists an integer $n\geq 0$, such that $K_n=K_{n+1}$, and what is the lower bound of $n$?

(3) Under what conditions do we know that the accumulation points come from dimension $1$ less. That is, can we find suitable $I',J'$, such that
\begin{center}
$\partial K(I,J,d)\subset K(I',J',d-1)$,
\end{center}
and when do we have
\begin{center}
$\partial K(I,J,d)=K(I,J,d-1)$?
\end{center}
\end{question}

For inductive reasons, we are interested in (1); for technical reasons, we are also interested in (2). Although (3) is stronger than both (1) and (2), (3) actually holds under some conditions. For example, we have the following for $LCT_{d}(I,J)$:
\begin{thm}\label{0104}\rm (Accumulation point theorem, \hyperref[HMX14]{[HMX14, Theorem 1.11]})
Suppose that $1$ is the only accumulation point of $I\subset [0,1]$, $I=I_+$ and $J=\mathbb N^+$. Then allaccumulation points of $LCT_d(I,J)$ are exactly $LCT_{d-1}(I)-\{1\}$. In particular, if $I\subset\mathbb Q$, the accumulation points of $LCT_{d}(I,J)$ are rational numbers.
\end{thm}

 In this paper we prove a similar statement for generalized pairs.

\begin{thm}\label{0105}\rm\textbf{Main Theorem} (Accumulation point theorem for generalized log canonical thresholds)
Let $d>0$ be an integer. Suppose that $I, J\subset [0,\infty)$ are two DCC sets. Then there are two DCC sets $I'\supset I$ and $J'\supset J$ that only depends on $I$, $J$ and $d$, such that
\begin{center}
$\partial GLCT_d(I,J)\subset GLCT_{d-1}(I',J')$.
\end{center}
In particular,

(1) If $I\subset\mathbb Q$, $J\subset\mathbb Q$ and all the accumulation points of $I, J$ are rational numbers, then $I', J'\subset\mathbb Q$, all the accumulation points of $I'$ are rational numbers, and all the accumulation points of $GLCT_{d}(I,J)$ are rational numbers.

(2) If $I$ only has finitely many accumulation points and $J$ has no accumulation point except $\infty$, then there exists an integer $m>0$ (that may depends on $I,J,d$) such that we can only take $m$ times accumulation points of $GLCT_{d}(I,J)$, i.e. there is no $(m+1)$-th order accumulation point of $GLCT_{d}(I,J)$.
\end{thm}

Indeed, we may prove a stronger statement than \hyperref[0105]{Theorem 1.5}:

\begin{thm}\label{0106}\rm
Let $d>0$ be an integer. Suppose that $I_1, I_2, J\subset [0,\infty)$ are three DCC sets. Let $GLCT_{d}(I_1,I_2,J)$ be the set
\begin{center}
$\{glct(X',B'+M';C'+N')|(X',B'+M')$ is generalized LC of dimension $d$,

$B'\in I_1, C'\in J, M=\sum m_jM_j, N=\sum n_jN_j$ where each $M_j,N_j$ is nef Cartier, $ m_j\in I_2, n_j\in J\}$
\end{center}

Suppose

(i) $1\in I_1=(I_1)_+$;

(ii) $I_2\cap [0,1]\subset I_1$;

(iii) $J$ is closed under addition.

Then there exists two DCC sets $I_1',J'$ and an integer $\mathcal{F}>0$, that only depends on $I_1,I_2,J, d$, such that

(1) The accumulation points of $GLCT_d(I_1,I_2; J)$ is either contained in $GLCT_{d-1}(\bar I_1,I_2; \bar J)$, or contained in $GLCT_{1}(I_1',\{0\};J')$;

(2) $J'=\frac{1}{\mathcal{F}}\bar J=\{\frac{j}{\mathcal{F}}|j\in \bar J\}$, $I_1'=(\frac{1}{\mathcal{F}}(\bar I_1\cup\bar I_2))_+$.

\end{thm}

\begin{thm}\label{0107}\rm
Let $d>0$ be an integer, $I_1,I_2,J\subset [0,\infty)$ be three DCC sets. 
Suppose 

(i) $I_1=(I_1)_+\subset [0,1]$ such that the only accumulation point of $I_1$ is $1$; 

(ii) $I_2\cap [0,1]\subset I_1$, and $I_2$ does not have any accumulation point except $\infty$.

(iii) $J$ is closed under addition, and $J$ has no accumulation point except $\infty$.

Then
\begin{center}
$\partial GLCT_{d}(I_1,I_2;J)=GLCT_{d-1}(I_1,I_2;J)$.
\end{center}

In particular, 

(1) If $I_1, I_2, J\subset\mathbb Q$, then all the accumulation points of $GLCT_{d}(I_1,I_2,J)$ are rational numbers.

(2) The only $d$-th order accumulation point of $GLCT_{d}(I_1,I_2,J)$ is $0$, and there's no $(d+1)$-th order accumulation point.

(3) In particular, the theorem holds when $I_2=\{0\}$, $J=\mathbb N$.
\end{thm}

\begin{rem}\label{0108}\rm
Before the start of the proof, it is worth to explain why we do not pick arbitrarily DCC sets, but pick DCC sets with good properties for \hyperref[0106]{Theorem 1.6} and \hyperref[0107]{Theorem 1.7}. Indeed, there are several reasons: DCC set are introduced when we do adjunction. Indeed, even for pairs when the coefficient set of divisors are only contained in $\{0,1\}$, after doing adjunction, we always get a coefficient set of form $\{1-\frac{1}{m}|m\in\mathbb N^+\}$. However, for DCC sets of form $D(I)\cup\{1\}$ (cf. \hyperref[0206]{Definition 2.6}), it remains the same set after doing adjunction, and that is the why for \hyperref[0106]{Theorem 1.6}, we need the assumption that $1\in I=I_+$.

For generalized pairs, however, based on similar adjunction formulas (cf. \hyperref[0306]{Theorem 3.6} and \hyperref[0307]{Theorem 3.7}), in order to have a good property when doing adjunction, it is necessary to at least suppose the following:

(i) $1\in I_1=(I_1)_+$; (ii) $I_2\cap [0,1]\subset I_1$; (iii) $J$ is closed under addition.

These are exactly the assumptions of \hyperref[0106]{Theorem 1.6}.

However, under this ``wild" assumption, it is only possible to show that we can find certain $I_1',I_2'$ and $J'$, such that 
\begin{center}
$\partial GLCT_{d}(I_1',I_2';J)\subset GLCT_{d-1}(I_1',I_2';J')$.
\end{center}

Now a natural question is, when are $I_1,I_2,J$ proper enough, such that
\begin{center}
$\partial GLCT_{d}(I_1,I_2;J)\subset GLCT_{d-1}(I_1,I_2;J)$?
\end{center}

In order to do this, we need to control the singularities in codimension $1$, and try to make it is as good as possible. It is not possible to say that $I_1,I_2,J$ all have no accumulation points, since $1$ is an accumulation point of $D(\{0,1\})$. However, notice that under adjunction, $I_2$ and $J$ remain ``unchanged" for lower-dimensional coefficient sets, and $I_1$ will ``change" to $D(I)$, we can assume that

(iv) $1$ is the only accumulation point of $I_1$, and 

(v) $I_2$, $J$ does not have any accumulation point (except $\infty$).

Together with (i), (ii), (iii), these are exactly the assumptions of \hyperref[0107]{Theorem 1.7}.

\end{rem}

\begin{flushleft} \textbf{Acknowledgement}. The author would like to thank his advisor Christopher D. Hacon for his constant support and many useful discussions. The author was partially supported by NSF research grants no: DMS-1300750, DMS-1265285 and by a grant from the Simons Foundation; Award Number: 256202.\end{flushleft}

\section{Notions And Conventions}

We will always work over the field of complex numbers $\mathbb C$.

\begin{defn}\rm (Positivity definitions) Let $X$ be a normal variety. An $\mathbb R$-Cartier $\mathbb R$-divisor $D$ on $X$ is called \textit{nef} if $D\cdot C\geq 0$ for any $C\subset X$; two $\mathbb R$-Cartier $\mathbb R$-divisors (resp. $\mathbb Q$-Cartier $\mathbb Q$-divisors) $D_1$, $D_2$ on $X$ are called $\mathbb R$ (resp.$\mathbb Q$)-\textit{linear equivalent} if $D_1-D_2$ is an $\mathbb R$ ($\mathbb Q$)-linear combination of principal divisors on $X$, and we use the notation $D_1\equiv_{lin,\mathbb R(\mathbb Q)}D_2$. Two $\mathbb R$-Cartier $\mathbb R$-divisors $D_1$, $D_2$ are called \textit{numerical equivalent} if $D_1\cdot C=D_2\cdot C$ for any curve $C\subset X$.

An $\mathbb R$-divisor $D$ on $X$ is called \textit{big} if $D\equiv_{lin,\mathbb R}A+C$ for some ample $\mathbb R$-divisor $A$ and some effective $\mathbb R$-divisor $C$. An $\mathbb R$-divisor $D$ is called \textit{pseudoeffective} if $D$ is a limit of big $\mathbb R$-divisors. An $\mathbb R$-divisor $D$ on $X$ is called \textit{nef} if $D\cdot C\geq 0$ for any $C\subset X$.

For any birational morphism $f: Y\rightarrow X$, an $\mathbb R$-divisor $D$ on $Y$ is called \textit{nef over $X$}, or $f$-\textit{nef}, if $D\cdot C\geq 0$ for any curve $C$ contained in $Ex(f)$. $D$ is called \textit{big over $X$}, or $f$-\textit{big}, if $D\equiv_{lin,\mathbb R}A+C$ for some $\mathbb R$-divisor $A$ that is ample over $X$, and some effective $\mathbb R$-divisor $C$.

An ample $\mathbb R$ (resp. $\mathbb Q$)-divisor $A$ is called \textit{general} if $A$ is a general element of $|A|_{\mathbb R}$ (resp. $|A|_{\mathbb Q}$), and \textit{very general} if $A$ is a very general element of $|A|_{\mathbb R}$ (resp. $|A|_{\mathbb Q}$).

A birational map $f: X\dashrightarrow Y$ is called $D$-nonpositive (resp. $D$-negative) for some $\mathbb R$-divisor $D$ on $X$ if $f$ does not extract any divisors, and for any common resolution of indeterminacy $p:W\rightarrow X$ and $q:W\rightarrow Y$, we have $p^*D=q^*D_Y+E$, such that $E\geq 0$ is $q$-exceptional (resp. $E\geq 0$ is $q$-exceptional and Supp$(p_*E)$ contains all the $f$-exceptional divisors), where $D_Y$ is the birational transform of $D$ on $Y$.
\end{defn}

\begin{defn}\label{0202}\rm (Pairs)
 A \textit{pair} $(X,B)$ is a normal variety $X$ and an effective $\mathbb R$-divisor $B$ on $X$. 
 
Let $\nu$ be a valuation of $X$ and $E$ be a prime divisor over $X$. $\nu$ (resp. $E$) is called exceptional over $X$ if $center_X\nu$ (resp. $center_XE$) is not a divisor. For any valuation $\nu$ of $X$ or any prime divisor $E$ over $X$, let $g: Y\rightarrow X$ be a log resolution of $(X,B)$ such that $center_Y\nu$ is a divisor (resp. $E$ is a divisor on $Y$). We define the \textit{log discrepancy}  $a(\nu, X,B)$ to be the number such that

\begin{center}
$K_Y=g^*(K_X+B)+\sum_{E_i=center_Y\nu_i\ \rm{is\ a\ divisor}} (a(\nu_i,X,B)-1)E_i$,
\end{center}
and if $E$ corresponds to a valuation $\nu$, we define $a(E,X,B)=a(\nu,X,B)$. 

 $(X,B)$ is called \textit{kawamata log terminal}, or \textit{KLT}, if $a(\nu,X,B)>0$ for any valuation $\nu$ of $X$. 
 
 $(X,B)$ is called \textit{divisorially log terminal}, or \textit{DLT}, if $a(\nu,X,B)>0$ for any valuation $\nu$ of $X$ such that $center_X\nu$ is not contained in the log smooth locus of $(X,B)$. 

 $(X,B)$ is called \textit{log canonical}, or \textit{LC} if $a(\nu,x,B)\geq 0$ for any valuation $\nu$ of $X$.

\end{defn}

\begin{defn}\label{0203}\rm (Generalized Pairs)
A \textit{generalized pair with data} $(M, X\xrightarrow{f} X'\xrightarrow{\pi} Z)$, or if without confusion, a \textit{generalized pair}, is a pair $(X',B')$ and a $\mathbb R$-divisor $M'$ written in the form $(X',B'+M')$, such that there exists projective morphisms $f: X\rightarrow X'$ and $\pi: X'\rightarrow Z$, a divisor $M$ on $X$, such that

(1) $f_*M=M'$;

(2) $f$ is birational;

(3) $X$ is normal;

(4) $M$ is nef over $Z$;

(5) $B'$ is a boundary, $M'$ is $\mathbb R$-Cartier.

$B'$ is called the \textit{boundary part} of $(X',B'+M')$ over $Z$, $M'$ is called the \textit{nef part} of $(X',B'+M')$ over $Z$, $Z$ is called the \textit{base} of the generalized pair $(X',B'+M')$, $(M, X\xrightarrow{f} X'\xrightarrow{\pi} Z)$ is called the \textit{data} of $(X',B'+M')$, $M$ is called the \textbf{nef data} of $(X',B'+M')$, and we say $(X',B'+M')$ is a \textit{generalized pair over $Z$}.

Sometimes we also say that $(X,B+M)$ is the data of $(X',B'+M')$.

If $f$ is a log resolution of $(X',B')$, we say the data is \textit{resolved}. If $Z=\{pt\}$ we possibly drop $Z$ and say that the generalized pair is \textit{projective}. If $Z=X'$, or $Z$ is not important, we also possibly drop $Z$. We may also disregard the data of the generalized pair if it is clear or not important.

\end{defn}

\begin{defn}\label{0204}\rm (Singularities of generalized pairs)
Let $(X',B'+M')$ be a generalized pair over $Z$, and possibly assuming it has a resolved data $(X,B+M)$. $(X',B'+M')$ is called \textit{generalized KLT} if all the coefficients of $B$ are $<1$, generalized log canonical or generalized LC if all the coefficients of $B$ are $\leq 1$, \textit{generalized} $\epsilon$-\textit{LC} if all the coefficients of $B$ are $\leq 1-\epsilon$. For any valuation $\nu$ on $X'$, possibly replacing $X$ by a resolved data such that $E=center_{X}\nu$ be a divisor, the \textit{generalized log discrepancy} $ga(\nu,X',B'+M')$ is defined as $1-\mu_EB$. We also use the notation $ga(E,X',B'+M')$. 

The \textit{generalized minimal log discrepancy} of $(X',B'+M')$, is defined as
\begin{center}
$min\{ga(\nu,X',B'+M')|\nu$ is exceptional over $X'\}$
\end{center}
and is denoted by $gmld(X',B'+M')$.

A \textit{generalized log canonical place} of a generalized log canonical pair $(X',B'+M')$ is a valuation $\nu$ on $X'$ such that $ga(\nu,X',B'+M')=0$. $center_{X'}\nu$ is called a \textit{generalized log canonical center}. $(X',B'+M')$ is called \textit{generalized DLT}, if $(X',B'+M')$ is generalized LC, $(X',B')$ is DLT, and every generalized log canonical center is a log canonical center of $(X',B')$. 

A \textit{generalized DLT modification} of a generalized log canonical pair $(X',B'+M')$ is a generalized pair $(X'',B''+M'')$ such that

(i) Possibly passing to a common resolved data, $(X',B'+M')$ and $(X'',B''+M'')$ has the same nef data $M$ and the same base $Z$.

(ii) there exists a birational morphism $\psi: X''\rightarrow X'$ such that
\begin{center}
$K_{X''}+B''+M''=\psi^*(K_{X'}+B'+M')$;
\end{center}

(iii) $(X'',B''+M'')$ is $\mathbb Q$-factorial generalized DLT;

(iv) $\psi$ only extracts valuations $\nu$ on $X'$ such that $ga(\nu,X',B'+M')=0$.
\end{defn}

\begin{defn}\label{0205}\rm (Adjunction of generalized pairs)
Let $(X',B'+M')$ be a generalized pair with data $(M,X\xrightarrow{f} X'\rightarrow Z)$. Let $S'\subset X'$ be a component of $\lfloor B'\rfloor$. We defined a generalized pair $(S',B_{S'}+M_{S'})$ in the following way:

Possibly after replacing $X$, we may assume the data is resolved, and we have
\begin{center}
$K_X+B+M=f^*(K_{X'}+B'+M')$.
\end{center} 

Let $S$ be the strict transform of $S'$ on $X$, $f|_S$ be the induced morphism $S\rightarrow S'$. We define $M_{S}=M|_S$, $B_S=(B-S)|_S$, $M_{S'}=(f|_{S})_*M_S$ and $B_{S'}=(f|_{S})_*B_S$. Thus, $(S',B_{S'}+M_{S'})$ is a generalized pair with data $(M_S,S\xrightarrow{f|_S} S'\rightarrow Z)$.
\end{defn}

\begin{defn}\label{0206}\rm (DCC and ACC sets)
Let $I$ be a set of real numbers. $I$ is called a \textit{DCC set} if any decreasing sequence contained in $I$ terminates. $I$ is called an \textit{ACC set} if any increasing sequence contained in $I$ terminates. $s$ is called an \textit{accumulation point of $I$ (from below (resp. above))} if  there exists a (strict increasing (resp. strict decreasing)) sequence contained in $I$ that converges to $s$. 

We define $\partial I$ to be the set 
\begin{center}
$\{s\in\mathbb R|s$ is an accumulation point of $I\}$
\end{center}
and $\bar I$ to be the set $I\cup\partial I$.

Let $m\geq 0$ be an integer. A real number $s$ is called the \textit{$m$-th order accumulation point of $I$} if $s\in\partial^mI$, where we define $\partial^0I=I$ and $\partial^mI=\partial\partial^{m-1}I$ inductively.

For any set $I\subset [0,\infty)$, we define the set $I_+$ as
\begin{center}
$\{0\}\cup\{j\in [0,1]|j=\sum_{p=1}^{l}i_p$ for some $l\in\mathbb N^+$ and $i_1,\dots,i_l\in I\}$
\end{center}
and the set $D(I)$ as
\begin{center}
$\{a\leq 1|a=\frac{m-1+f}{m}, m\in\mathbb N^+,f\in I_+\}$.
\end{center}

For any set $I\subset [0,1]$, we define the set $\Phi(I)$ as
\begin{center}
$\{1-\frac{r}{m}|r\in I, m\in\mathbb N^+\}$.
\end{center}

For any set $I\subset [0,\infty)$, any real number $c\geq 0$, we define the set $c\cdot I$, or $cI$, as
\begin{center}
$\{ci|i\in I\}$.
\end{center}

A set $I\subset [0,\infty)$ is called \textit{closed under addition} if for any $i,j\in I$, $i+j\in I$. For example, $\mathbb N$ is closed under addition.

Let $X'$ be a normal variety, $B$ an $\mathbb R$-divisor on $X'$. We write $B\in I$ if all the coefficients of $B$ are contained in $I$. 

Let $f: X\rightarrow X'$ be a birational morphism, $\pi: X'\rightarrow Z$ be a contraction such that $f$ and $\pi$ are both projective. Let $M$ be a divisor on $X$ such that $M$ is nef over $Z$. We write $M\in\mathcal{NEF}_{X/Z}(I)$ if we can write $M=\sum m_jM_j$, such that each $M_j$ is a nef Cartier divisor over $Z$, each $m_j\in I$. If $Z$ is a point or not important, we may disregard $Z$ and write $M\in\mathcal{NEF}_{X}(I)$.
\end{defn}

\begin{defn}\label{0207}\rm (Thresholds)

Suppose $(X',B'+M')$ is a generalized log canonical pair over $Z$ with data $(M,X\xrightarrow{f} X'\rightarrow Z)$. Let $N'$ be a $\mathbb R$-Cartier divisor on $X'$, such that after possibly replacing $X$ by a resolved data, there exists a nef $\mathbb R$-divisor $N$ on $X$ such that $f_*N=N'$. Let $C'$ be an effective $\mathbb R$-divisor on $X'$. The \textit{generalized log canonical threshold} of $(X',B'+M')$ respect to $C'+N'$ is defined as
\begin{center}
$\underset{c\geq 0}{sup}\{(X',B'+cC'+M'+cN')$ is generalized log canonical with nef data $M+cN\}$
\end{center}
and is denoted by
\begin{center}
$glct(X',B'+M';C'+N')$.
\end{center}

If $M=N=0$, $X=X'=Z$, we define
\begin{center}
$lct(X',B';C')=glct(X',B'+M',C'+N')$.
\end{center}

Let $I,J\subset [0,\infty)$ be two sets, $d>0$ be an integer. We define
\begin{center}
$LCT_d(I;J):=\{lct(X,B;C)|(X,B)$ is log canonical of dimension $d, B\in I, C\in J\}$ 
\end{center}

\begin{center}
$GLCT_{d}(I;J):=\{glct(X',B'+M';C'+N')|(X',B'+M')$ is generalized LC of dimension $d$,

$(M,X\rightarrow X')$ is the data, $B'\in I, C'\in J, M\in\mathcal{NEF}_{X/X'}(I), N\in\mathcal{NEF}_{X/X'}(J)\}$.
\end{center}

Let $I_1,I_2,J\subset [0,\infty)$ be three sets, $d>0$ be an integer. We define
\begin{center}
$GLCT_{d}(I_1,I_2;J):=\{glct(X',B'+M';C'+N')|(X',B'+M')$ is generalized LC of dimension $d$,

$(M,X\rightarrow X')$ is the data, $B'\in I_1, C'\in J, M\in\mathcal{NEF}_{X/X'}(I_2), N\in\mathcal{NEF}_{X/X'}(J)\}$.
\end{center}
\end{defn}

\begin{defn}\label{0208}\rm (Definitions for technical requirements)

Let $I_1,I_2,J\subset [0,\infty)$ be three sets. Let $d>0$ be an integer, $c>0$ be a real number. We define
\begin{center}
$\mathcal{D}_c(I_1; J)=\{a\leq 1|a=\frac{m-1+f+kc}{m}, m\in\mathbb N^+, f\in (I_1)_+, k=\sum_{p=1}^{l} j_p\not=0$ for some $l$ and $j_1,\dots, j_l\in J\}$;
\end{center}

\begin{center}
$\mathcal{R}_d(I_1,I_2;J,c)=\{(X',\Delta'+\Gamma')|(X',\Delta'+\Gamma')$ is projective generalized log canonical of dimension $d$,

 with boundary part $\Delta'$ and nef data $\Gamma$, $K_{X'}+\Delta'+\Gamma'\equiv_{num}0$,
 
 $\Delta'=B'+C'$ where $B'\in D(I_1)$ and $C'\in\mathcal{D}_c(I_1;J)$;
 
 $\Gamma=M+cN$ where $M\in\mathcal{NEF}_{X}(I_2)$ and $N\in\mathcal{NEF}_{X}(J),$
 
 such that either $N\not\equiv_{num}0$, or $C'\not=0\}$;
\end{center}

\begin{center}
$\mathcal{S}_d(I_1,I_2;J,c)=\{(X',\Delta'+\Gamma')\in\mathcal{R}_n(I_1,I_2;J,c)$ for some $n\leq d, $

$\rho(X')=1, (X',\Delta'+\Gamma')$ is $\mathbb Q$-factorial generalized KLT$\}$;
\end{center}

\begin{center}
$\mathcal{N}_d(I_1,I_2;J)=\{c\geq 0|\mathcal{R}_d(I_1,I_2;J,c)\not=\emptyset\}$;
\end{center}
and

\begin{center}
$\mathcal{K}_d(I_1,I_2;J)=\{c\geq 0|\mathcal{S}_d(I_1,I_2;J,c)\not=\emptyset\}$.
\end{center}
\end{defn}

If $I_1=I_2$, we may replace each of the above notations by $\mathcal{D}_c(I;J), \mathcal{R}_d(I;J,c),\mathcal{S}_d(I;J,c),\mathcal{N}_d(I;J)$ and $\mathcal{K}_d(I;J)$.

\section{Preliminaries}

In this section we state the theorems that we will repeatedly use in our proof of the main theorem. First we shall state the ACC theorems for pairs and generalized pairs and their global version.
\begin{thm}\label{0301}\rm (ACC for log canonical thresholds)
Let $d>0$ be an integer, $I\subset [0,1]$, $J\subset [0,\infty)$ be two DCC sets. Then $LCT_d(I;J)$ satisfies the ACC.
\end{thm}

\begin{flushleft}\textit{Proof}. \hyperref[HMX14]{[HMX14,Theorem 1.1]}.\hfill$\Box$\end{flushleft}

\begin{thm}\label{0302}\rm (ACC for log canonical thresholds, global version)
Let $d>0$ be an integer, $I\subset [0,1]$ be a DCC set. Then there exists a finite subset $I_0\subset I$ that only depends on $I$ and $d$, such that for any projective log canonical pair  $(X,B)$ of dimension $d$, if $K_X+B\equiv_{num}0$ and $B\in I$, then $B\in I_0$.
\end{thm}

\begin{flushleft}\textit{Proof}. \hyperref[HMX14]{[HMX14,Theorem 1.5]}.\hfill$\Box$\end{flushleft}

\begin{thm}\label{0303}\rm (ACC for generalized log canonical thresholds)
Let $d>0$ be an integer, $I, J\subset [0,\infty)$ be two DCC sets. Then $GLCT_{d}(I,J)$ satisfies the ACC.
\end{thm}

\begin{flushleft}\textit{Proof}. \hyperref[BZ16]{[BZ16, Theorem 1.5]}.\hfill$\Box$\end{flushleft}

\begin{thm}\label{0304}\rm (ACC for generalized log canonical thresholds, global version)
Let $d>0$ be an integer, $I\subset [0,\infty)$ be a DCC set. Then there exists a finite subset $I_0\subset I$ that only depends on $I$ and $d$, such that for any projective generalized log canonical pair $(X',B'+M')$ of dimension $d$ with nef part $M$, if 

(i) $M=\sum m_jM_j$, such that each $M_j$ is nef Cartier, each $m_j\in I$, and if $m_j\not=0$, $M_j\not\equiv_{num}0$;

(ii) $B'\in I$;

(iii) $K_{X'}+B'+M'\equiv_{num}0$.

Then $B'\in I_0$, $M\in\mathcal{NEF}_{X}(I_0)$.
\end{thm}

\begin{flushleft}\textit{Proof}. \hyperref[BZ16]{[BZ16, Theorem 1.6]}.\hfill$\Box$\end{flushleft}

Next we shall state several results of generalized adjunction, which are similar to \hyperref[BZ14]{[BZ14, 4.7-4.10]}.

\begin{thm}\label{0305}\rm (Generalized log canonicity under adjunction)
Let $(X',B'+M')$ be a generalized log canonical pair, $S'\subset X'$ be a component of $\lfloor B'\rfloor$. Then $(S',B_{S'}+M_{S'})$ given by the generalized adjunction is generalized log canonical, and in particular, $B_{S'}$ is a boundary. 
\end{thm}

\begin{flushleft}\textit{Proof}. \hyperref[BZ14]{[BZ14, 4.7-4.8]}.\hfill$\Box$\end{flushleft}

\begin{thm}\label{0306}\rm (Coefficient of generalized adjunction)
Let $(X',B'+M')$ be a generalized log canonical pair with data $(M, X\xrightarrow{f}X'\rightarrow Z)$, $1\in I\subset [0,\infty)$ be a set. Let $S'\subset X'$ be a component of $\lfloor B'\rfloor$. Suppose $B'\in I$ and $M\in\mathcal{NEF}_{X}(I)$. Then all the coefficients of $B_{S'}$ given by the generalized adjunction
\begin{center}
$K_{S'}+B_{S'}+M_{S'}=(K_{X'}+B'+M')|_{S'}$
\end{center}
are contained in the set $D(I)$, and $M_{S'}\in\mathcal{NEF}_{S}(I)$.
\end{thm}

\begin{flushleft}\textit{Proof}. To show $B_{S'}\in D(I)$, a similar but non-precise argument can be found in \hyperref[BZ14]{[BZ14, Proposition 4.9]}. We shall give a direct proof.\end{flushleft}

Possibly replacing $X$ we may assume the data is resolved. Let $S$ be the strict transform of $S'$ on $X$, and let $V'$ be a component of $B_{S'}$. We consider the image of the generic point $\eta_{V'}$ of $V'$ in $X'$. By \hyperref[Sho93]{[Sho93, Corollary 3.9]}, there exists an integer $m>0$, such that for any Weil divisor $D'$ on $X'$, $mD'$ is Cartier near $\eta_{V'}$. Let 
\begin{center}
$K_{S'}+\tilde B_{S'}=(K_{X'}+B')|_{S'}$
\end{center}
be the usual adjunction. Then by \hyperref[Sho93]{[Sho93, Corollary 3.10]}, the coefficient of $V'$ in $\tilde B_{S'}$ is of the form
\begin{center}
$\frac{m-1+i}{m}$
\end{center}
where $i\in I_+$. Moreover, suppose $M=\sum m_jM_j$ where each $M_j$ is nef Cartier over $Z$ and $m_j\in I$, and let $M_j'=f_*M_j$ for every $j$. Then near the generic point of $V'$, for any $j$, we have
\begin{center}
$f^*M_j'=M_j+E_j$
\end{center}
for some exceptional divisor $E_j$ over $X'$. Since $M_j$ is nef, by the negativity lemma, $E_j$ is effective; since $M_j$ is Cartier, near the generic point of $V'$, $mE_j$ for every $j$. Thus for every $j$, the multiplicity of $E_j|_{S_i}$ near the image of the generic point of $V'$ is of the form
\begin{center}
$\frac{z_j}{m}$
\end{center}
where $z_j\in\mathbb N$. Thus, since
\begin{center}
$M_{S'}+(B_{S'}-\tilde B_{S'})=M'|_{S'}$,
\end{center}
the coefficient of $V'$ contained in $B_{S'}$ is of the form
\begin{center}
$\frac{m-1+i+\sum m_jz_j}{m}\in D(I)$.
\end{center}

$M_{S'}\in\mathcal{NEF}_{S}(I)$ is straightforward from definition.\hfill$\Box$

\begin{thm}\label{0307}\rm (Coefficient of generalized adjunction with particular coefficient)

Let $I,J\subset [0,\infty)$ be two sets. $(X',\Delta'+\Gamma')$ be a generalized pair with data $(\Gamma=M+N, X\xrightarrow{f}X'\rightarrow Z)$ such that $\Delta'=B'+C'$ is the boundary part and $\Gamma'=M'+N'$ is the nef part, where $M'=f_*M$ and $N'=f_*N$. Suppose $B'\in I$ and $M\in\mathcal{NEF}_{X/Z}(I)$, $C'\in J$ and $N\in\mathcal{NEF}_{X/Z}(J)$. 

Let $S'$ be a component of $\lfloor B'\rfloor$. Let
\begin{center}
$K_{S'}+\Delta_{S'}+\Gamma_{S'}=(K_{X'}+\Delta'+\Gamma')|_{S'}$
\end{center}
be the generalized adjunction.

Suppose either $C'|_{S'}\not=0$ or $N'|_{S'}\not\equiv_{num,Z}0$. Then:

1) Either there exists a nonzero component of $\Delta_{S'}$ of form
\begin{center}
$\frac{m-1+i+\sum_{p=1}^{l} j_pz_p}{m}$
\end{center}
where $i\in I_+$, $l\geq 1$, $j_p\in J$ for each $p$, $z_p\in\mathbb N^+$ for each $p$, $m\in\mathbb N^+$, or 

2) $N|_{S}\not\equiv_{num,Z}0$. In particular, $N_{S}=\sum n_kN_k|_{S}$ where $n_k\in J\backslash\{0\}$, and there exists $k$ such that $N_{k}|_{S}\not\equiv_{num,Z}0$.
\end{thm}

\begin{flushleft}\textit{Proof}. This is similar to \hyperref[BZ14]{[BZ14, Lemma 4.10]} but we just state it in a more accurate way. By \hyperref[0306]{Theorem 3.6}, all the components of $\Delta_{S'}$ are of the form 
\begin{center}
$\frac{m-1+s}{m}$
\end{center}
where $s\in (I\cup J)_+$. Thus $s=i+j$, where $i\in I_+$ and $j\in J_+$. If $j\not=0$, we are in case 1) and we are done. Thus we may assume $j=0$, and in particular, all the coefficients of $\Delta_{S'}$ are contained in the set $D(I)$.\end{flushleft}

For any component $V'$ of $\Delta_{S'}$, we consider the image of the generic point $\eta_{V'}$ of $V'$ in $X'$. Suppose near $\eta_{V'}$, $mD'$ is Cartier for any Weil divisor $D'$ on $X'$. Let 
\begin{center}
$K_{S'}+\tilde \Delta_{S'}=(K_{X'}+\Delta')|_{S'}$
\end{center}
be the usual adjunction. Then by \hyperref[Sho93]{[Sho93, Corollary 3.10]}, the coefficient of $V'$ contained in $\tilde\Delta_{S'}$ is of the form 
\begin{center}
$\frac{m-1+i+j}{m}$
\end{center}
where $i\in I_+$, $j\in J_+$, and if $C'|_{S'}\not=0$, $j\not=0$.

Now suppose $M=\sum m_kM_k$ and $N=\sum n_kN_k$ where each $M_k,N_k$ is nef Cartier over $Z$, $m_k\in I$ and $n_k\in J$. Let $M'_k=f_*M_k$ and $N'_k=f_*N_k$. 

Suppose $f^*M'_k=M_k+E_k$ and $f^*N'_k=N_k+F_k$ where $E_k,F_k$ are exceptional divisors. By negativity lemma, $E_k,F_k$ are effective, and since $M_k,N_k$ are Cartier, near the image of the generic point of $V'$, $mE_k$ and $mF_k$ are Cartier. Thus the multiplicity of $V'$ of $E_k|_{S}$ and $F_k|_S$ are of the form $\frac{z_k}{m}$ and $\frac{z_k'}{m}$ where $z_k,z_k'\in\mathbb N$, and $z_k'\not=0$ if $F_k|_{S}\not=0$.

Now the coefficient of $V'$ contained in $\Delta_{S'}$ is of the form
\begin{center}
$\frac{m-1+i+j+\sum z_km_k+\sum z_k'n_k}{m}$
\end{center}
where $i\in I_+, j\in J_+, z_k,z_k'\in\mathbb N$, $m_k\in I, n_k\in J$. Since $\frac{m-1+i+j+\sum z_km_k+\sum z_k'n_k}{m}\in D(I)$, $j=0$ and each $z_k'=0$. In particular, $C'|_{S'}=0$, thus $N'|_{S'}\not\equiv_{num,Z}0$. But since $F_k|_{S}=0$, $N|_{S}\not\equiv_{num}0$, and we are in case 2).\hfill$\Box$\vspace{2mm}

We state the theorem that generalized divisorial log terminal modification exists, which is needed in many of our following statements.

\begin{thm}\label{0308}\rm

Let $(X',B'+M')$ be a generalized LC pair. Then there exists a generalized DLT modification of $(X',B'+M')$. In particular, if $(X',B'+M')$ is generalized DLT, it is generalized KLT if and only if $\lfloor B'\rfloor=0$.
\end{thm}

\begin{flushleft}\textit{Proof}. The existence of generalized DLT modification is guaranteed by \hyperref[BZ14]{[BZ14, Lemma 4.5]}. Suppose $(X',B'+M')$ is generalized DLT. If it is generalized KLT, $\lfloor B'\rfloor=0$ by definition. If $\lfloor B'\rfloor=0$, then since $(X',B')$ is DLT, $(X',B')$ is KLT, hence there does not exist a generalized LC center of $(X',B'+M')$, hence $(X',B'+M')$ is generalized KLT.\hfill$\Box$\end{flushleft}

We show a boundedness theorem on the number of components.

\begin{thm}\label{0309}\rm (Bound of number of components)
Let $d>0$ be an integer, $b,c\geq 0$ be two positive numbers. Then there exists an integer $p>0$ that only depends on $d$ and $b$, such that for any projective generalized pair $(X',B+M')$ with data $(M,X\xrightarrow{f}X')$, such that 

(i) $M=\sum_{j=1}^r m_jM_j$, where each $M_j$ is Cartier;

(ii) Each $M_j':=f_*M_j$ is big;

(iii) $K_{X'}+B'+M'+P'\equiv_{num}0$ for some pseudoeffective $\mathbb R$-divisor $P'$;

(iv) All the coefficients of $B'$ are $\geq c$.

Then the number of components of $B'$ is at most $\frac{d+1}{c}$ and $r\leq p$.
\end{thm}

\begin{flushleft}\textit{Proof}. \hyperref[BZ14]{[BZ14, Proposition 5.1 and 5.2]}.\hfill$\Box$\end{flushleft}

We show an easy lemma that is helpful in many instances

\begin{lem}\label{0310}\rm
Let $X'\rightarrow Z$ be a projective morphism. Let $X\xrightarrow{f} X'$ be a birational morphism, $N$ be a nef $\mathbb R$-Cartier $\mathbb R$-divisor over $Z$ on $X$, and $N'$ be the pushforward of $N$ to $X'$. Suppose that $N'$ is $\mathbb R$-Cartier. Then:

(1) If $N$ is numerically trivial over $Z$, so is $N'$.

(2)  If $Z=\{pt\}$, then $N$ is numerically trivial iff $N'$ is numerically trivial.
\end{lem}

\begin{flushleft}\textit{Proof}. Suppose $f^*N'=N+E$. Then $E$ is exceptional over $X'$, and by the negativity lemma, $E\geq 0$. To prove (1), since $N$ is numerically trivial, by the negativity lemma again, $E\leq 0$. Thus $E=0$, and $f^*N'=N$ is numerically trivial. Thus by the projection formula, $N'$ is numerically trivial over $Z$.\end{flushleft}

To prove (2), if $N'$ is numerically trivial, $N\equiv_{num}-E$. If $E=0$ then we are done. Otherwise, let $x'$ be a point that contained in the center of a component of $E$ on $X'$. Now we pick a general curve $\Sigma'$ contained in $X'$ that passes through $x'$, and let $\Sigma$ be the strict transform of $\Sigma'$ on $X$. Then $N\cdot\Sigma'=-E\cdot\Sigma'<0$, which is not possible.\hfill$\Box$\vspace{2mm}

Finally, we need the following boundedness result proved by Birkar:

\begin{thm}\label{0311}\rm
Let $d>0$ be an integer, $\epsilon>0$ be a real number. Then the set of all $\epsilon$-log canonical Fano varieties forms a bounded family.
\end{thm}

\begin{flushleft}\textit{Proof}. \hyperref[Bir16]{[Bir16, Theorem 1.1]}.\end{flushleft}

\section{Relationship Between Sets for Technical Requirements}

In this section we deal with the relationships between the sets defined in \hyperref[0208]{Definition 2.8}. We show similar statements as in \hyperref[HMX14]{[HMX14, Section 11]}.

\begin{lem}\label{0401}\rm
Let $I_1,I_2,J\subset [0,\infty)$ be three sets, $d>0$ be an integer, and $c>0$ be a real number. Then

(1) $GLCT_{d}(I_1,I_2;J)\subset GLCT_{d+1}(I_1,I_2;J)$;

(2) $\mathcal{N}_d(I_1,I_2;J)\subset\mathcal{N}_{d+1}(I_1,I_2;J)$;

(3) $\frac{1-s}{\sum j_kz_k}\in\mathcal{N}_d(I_1,I_2;J)$ and $GLCT_{d}(I_1,I_2;J)$ for any $s\in (I_1)_+$, $z_k\in\mathbb N$, $j_k\in J$ such that $\sum j_kz_k\not=0$;

(4) $\frac{1-i}{j}\in GLCT_d(I_1,I_2;J)$ for any $i\in I_1$ and $j\in J$.
\end{lem}

\begin{flushleft}\textit{Proof}. Let $E$ be an elliptic curve, $(X',B'+M')$ be a generalized pair with data $(M,X\xrightarrow{f}X'\rightarrow Z)$. Then there exists a natural contraction $Z\times E\rightarrow Z$ and natural morphisms $X\times E\rightarrow X'\times E\rightarrow Z\times E$. Since $E$ is  projective, the morphisms remains projective; since $E$ is smooth, $X\times E$ is smooth. \end{flushleft}

Write $B_i=\sum b_iB_i$, $B_i'=\sum b_iB_i'$ as in irreducible components where $B_i'=f_*B_i$ for each $i$, and $M_i=\sum m_iM_i$, $M_i'=\sum m_iM_i'$ where $M_i'=f_*M_i$ for each $i$ such that $M_i$ is nef Cartier over $Z$. Then we consider the divisors $B_i\times E\subset X\times E$, $B_i'\times E\subset X'\times E$, $M_i\times E\subset X\times E$ and $M_i'\times E\subset X'\times E$. For each $i$, it is clear that $M_i\times E$ is nef over $Z\times E$, and since $E$ is an elliptic curve, $M_i\times E$ is nef over $Z$. Thus, $(X'\times E, B'\times E+M'\times E:=\sum b_i(B_i'\times E)+\sum m_i(M_i'\times E))$ has a structure of generalized pair with data $(M=\sum m_i(M_i\times E),X\times E\rightarrow X'\times E\rightarrow Z)$. 

Since $E$ is an elliptic curve, $(X'\times E,B'\times E+M'\times E)$ is generalized log canonical (resp. generalized KLT) if and only if $(X',B'+M')$ is generalized log canonical (resp. generalized KLT), $K_{X'}+B'+M'\equiv_{num}0$ if and only if $K_{X'\times E}+B'\times E+M'\times E\equiv_{num}0$, and for any set $I\subset [0,\infty)$, $B'\in I$ iff $B'\times E\in I$, and $M\in\mathcal{NEF}_{X/Z}(I)$ iff $M\times E\in\mathcal{NEF}_{X\times E/Z}(I)$. In particular, (1) and (2) hold.

To prove (3), we only need to prove $\frac{1-f}{\sum j_kz_k}\in\mathcal{N}_1(I_1,\{0\};J)$, and then since $0\in I_2$, by using (2) we can show that $\frac{1-f}{\sum j_kz_k}\in\mathcal{N}_d(I_1,I_2;J)$.

Let $X'=X=\mathbb P^1, Z=\{pt\}$. For any integer $p>0$, any $j_1,\dots, j_p\in J\backslash\{0\}$ and $z_1,\dots, z_p\in\mathbb N^+$, any $s\in (I_1)+$, let $x_1,\dots,x_p, y_1,\dots, y_p, w,v$ be $2p+2$ different points on $X'$, we consider the generalized pair
\begin{center}
$(X',s(w)+s(v)+\sum_{l=1}^pc\cdot j_p(x_l)+\sum_{l=1}^pc\cdot j_p(y_l))$.
\end{center}
where $\Delta'=s(w)+s(v)$ is the boundary part and $\Gamma'=\sum_{l=1}^pc\cdot j_p(x_l)+\sum_{l=1}^pc\cdot j_p(y_l)$ is the nef part.

Since $\Delta'=B'+C'$ where $C'=0$ and $B'=\Delta'\in I_+\subset D(I)$, and $\Gamma'=0+N'$ where $N=N'\in\mathcal{NEF}_{X}(J)$ and $N\not\equiv_{num}0$, $\mathcal{R}_1(I_1,\{0\};J,c)\not=\emptyset$, hence $c\in\mathcal{N}_1(I_1,\{0\},J)$.

To prove (4), we only need to prove $c=\frac{1-i}{j}\in GLCT_1(I_1,\{0\},J)$. Let $X=X'=\mathbb P^1$, $M=N=0$, $p\in X$ be a general point, $B=ip$ and $C=jp$. Then
\begin{center}
$glct(X,B;C)=c$.
\end{center}
\hfill$\Box$

\begin{lem}\label{0402}\rm
Let $I_1,I_2,J\subset [0,\infty)$ be three sets, such that $1\in I_1$, and $I_2\cap [0,1]\subset I_1$. Then
\begin{center}
$\mathcal{N}_d(I_1,I_2;J)=\mathcal{K}_d(I_1,I_2;J)$.
\end{center}
\end{lem}

\begin{flushleft}\textit{Proof}. $\mathcal{N}_d(I_1,I_2;J)\supset\mathcal{K}_d(I_1,I_2;J)$ is clear from the definitions and \hyperref[0401]{Lemma 4.1(2)}. Thus, we only need to prove that $\mathcal{N}_d(I_1,I_2;J)\subset\mathcal{K}_d(I_1,I_2;J)$. \end{flushleft}

For every $c\in\mathcal{N}_d(I_1,I_2;J)$, let $(X',\Delta'+\Gamma')$ be a projective generalized log canonical pair with data $(\Gamma, X\rightarrow X')$ such that $\Delta'=A'+B'+C'$ where $A'$ is reduced (in particular, $A'\in D(I_1)$), $B'\in D(I_1)\backslash\{1\}$, and $C'\in D_c(I_1;J)$, $\Gamma=M+cN$ where $M\in\mathcal{NEF}_{X}(I)$ and $N\in\mathcal{NEF}_{X}(J)$, such that either $C'\not=0$ or $N\not\equiv_{num}0$, and $K_{X'}+\Delta'+\Gamma'\equiv_{num}0$.

By \hyperref[0308]{Theorem 3.8}, we may take a generalized DLT modification $(X'',\Delta''+\Gamma'')$ of $(X',\Delta'+\Gamma')$, and possibly replacing $X$ by a resolved data, we may assume that $X\dashrightarrow X''$ is a morphism. We may write
\begin{center}
$K_{X''}+A''+B''+C''+M''+cN''=\psi^*(K_{X'}+A'+B'+C'+M'+cN')$
\end{center}
such that $\Gamma''=M''+cN''$ is the nef part, $B''$, $C''$ are the strict transforms of $B'$ and $C'$, and $A''$ is the strict transform of $A'$ plus the reduced exceptional divisor of $\psi$. In particular, since $K_{X''}+\Delta''+\Gamma''\equiv_{num}0$, $(X'',\Delta''+\Gamma'')\in\mathcal{R}_d(I_1,I_2;J,c)$. Replace $(X',\Delta'+\Gamma')$ by $(X'',\Delta''+\Gamma'')$, we may assume that $(X',\Delta'+\Gamma')$ is $\mathbb Q$-factorial DLT.

Let $A$ be the strict transform of $A'$ on $X$. There are two cases:

\begin{flushleft}\textbf{Case 1}\label{1042}. $A'\not=0$, and there exists a component $S'$ of $A'$, such that either $C'$ intersects $S'$, or $N'|_{S'}\not\equiv_{num}0$. Notice that $(I_1\cup I_2)_+=(I_1)_+$, let 
\begin{center}
$K_{S'}+\Delta_{S'}+\Gamma_{S'}=(K_{X'}+\Delta'+\Gamma')|_{S'}$
\end{center}
be the generalized adjunction. By \hyperref[0306]{Theorem 3.6}, $\Gamma_{S}=M_S+cN_S$ where $M_S\in\mathcal{NEF}_{S}(I_2)$, $N_S\in\mathcal{NEF}_S(J)$, $\Delta_{S'}\in D(I_1\cup I_2\cup c\cdot J)=D(I_1\cup c\cdot J)$. Thus, $\Delta_{S'}$ can be written in the form of $B_{S'}+C_{S'}$ where $B_{S'}\in D(I_1)$, $C_{S'}\in D_c(I_1;J)$. \end{flushleft}

By \hyperref[0307]{Theorem 3.7}, either $C_{S'}\not=0$ or $N_S\not\equiv_{num}0$. Thus, $(S',\Delta_{S'}+\Gamma_{S'})$ is generalized LC of dimension $d-1$. By using induction on $d$, $c\in\mathcal{K}_{d-1}(I_1,I_2;J)\subset\mathcal{K}_d(I_1,I_2;J)$, and we are done.

\begin{flushleft}\textbf{Case 2}\label{2042}. No component $S'$ of $A'$ intersects $C'$ and $N'|_{S}\not=0$. Since $N$ is pseudoeffective, $N'$ is pseudoeffective. Thus if $C'\not=0$, $C'+cN'$ is pseudoeffective but not numerically trivial; and if $C'=0$, $N\not\equiv_{num}0$, hence $N'\not\equiv_{num}0$, and $C'+cN'$ is pseudoeffective but not numerically trivial.

In either case, $K_{X'}+A'+B'+M'\equiv_{num}-C'-cN'$ is not pseudoeffective, and we may run a $(K_{X'}+A'+B'+M')$-MMP. Let $\phi: X'\rightarrow X''$ be a step of a $(K_{X'}+A'+B'+M')$-MMP with scaling of some ample divisor. Then $\phi$ is $(C'+cN')$-positive.\end{flushleft}

\begin{flushleft}\textbf{Case 2.1}\label{2142} $\phi$ does not define a Mori fiber space. \end{flushleft}

First suppose we are in \hyperref[2142]{Case 2.1}. Possibly replacing $X$ by resolved data, we may assume $X\rightarrow X''$ is a morphism. Let $A'',B'',C'',M'',N''$ be the pushforward of $A,B,C,M,N$. Then $M\in\mathcal{NEF}_{X}(I_2)$, $N\in\mathcal{NEF}_{X}(J)$, $A''$ is reduced, $B''\in D(I_1)$, and $C''\in\mathcal{D}_c(I_1;J)$. Moreover, $K_{X''}+A''+B''+C''+M''+cN''\equiv_{num}0$. 

Notice that if $N\not\equiv_{num}0$, from our definition, $(X'',A''+B''+C''+M''+cN'')\in\mathcal{R}_{d}(I_1,I_2;J,c)$; and if $N\equiv_{num}0$, then $C''\not=0$, and $\phi$ is indeed a step of $(-C'')$-MMP, hence $C''$ is not contracted. In other words, $C''\not=0$, thus $(X'',A''+B''+C''+M''+cN'')\in\mathcal{R}_{d}(I_1,I_2;J,c)$.

If $\phi$ is a divisorial contraction or a flip, and if any component $S'$ of $A'$ is contracted, then either $C'$ intersects $A'$, or $N'|_{S'}\not\equiv_{num}0$. In this case we are back to \hyperref[1042]{Case 1} and the proof is finished. Otherwise, no component of $A'$ is contracted, we replace $(X',A'+B'+C'+M'+cN')$ by $(X'',A''+B''+C''+M''+cN'')$ and continue running a $(K_{X'}+A'+B'+M')$-MMP with scaling. If we are still in \hyperref[2142]{Case 2.1}, we repeat the process above, but we cannot repeat infinitely many times since this MMP must terminate with a Mori fiber space structure. Thus, repeating the process finitely many times, $\phi$ must defines a Mori fiber space, and we move on after \hyperref[2242]{Case 2.2}.

\begin{flushleft}\textbf{Case 2.2}\label{2242} $\phi$ defines a Mori fiber space.\end{flushleft}

If we are in \hyperref[2242]{Case 2.2}, suppose $\phi: X'\rightarrow T$ is the Mori fiber space. If $dim T>0$, let $t\in T$ be a general point, and $F'_t$ be a general fiber of $\phi$. Then adjunction to the fiber $F'=F'_t$ gives a generalized pair $(F',A_{F'}+B_{F'}+C_{F'}+M_{F'}+cN_{F'})\in\mathcal{R}_{d-dim T}(I_1,I_2;J,c)$, thus $c\in\mathcal{K}_{d-dim T}(I_1,I_2;J)\subset\mathcal{K}_{d}(I_1,I_2;J)$ and we are done. Hence we may assume $dim T=0$ and thus $\rho(X')=1$.

If $A'=0$, the pair $(X',A'+B'+C'+M'+cN')\in\mathcal{K}_d(I_1,I_2;J)$ by definition. Otherwise, $A'\not=0$. But since $X'$ is of Picard number $1$, $A'$ is ample hence either $C'$ intersects $A'$, or $N'|_{S'}\not\equiv_{num}0$ for any component $S'$ of $A'$, and we are back to \hyperref[1042]{Case 1}.\hfill$\Box$

\begin{lem}\label{0403}\rm
Let $d>0$ be an integer. Suppose $I_1,I_2,J\subset [0,\infty)$ are three sets, such that $1\in I_1=(I_1)_+$, $I_2\cap [0,1]\subset I_1$, and $J$ is closed under addition. Then 
\begin{center}
$GLCT_{d+1}(I_1,I_2;J)=\mathcal{N}_{d}(I_1,I_2;J)$.
\end{center}
\end{lem}

\begin{flushleft}\textit{Proof}. First we show that $GLCT_{d+1}(I_1,I_2;J)\subset\mathcal{N}_d(I_1,I_2;J)$. Let $c\in GLCT_{d+1}(I_1,I_2;J)$. Then there exists a generalized pair $(X',B'+M')$ with data $(M,X\rightarrow X'\rightarrow Z)$ such that $B'\in I_1$, $M\in\mathcal{NEF}_{X}(I_2)$, an $\mathbb R$-divisor $C'\in J$ contained in $X'$, a nef $\mathbb R$-divisor $N\in\mathcal{NEF}_{X}(J)$, such that
\begin{center}
$c=glct(X',B'+M';C'+N')$
\end{center}
where $N'$ is the pushforward of $N$ to $X'$. Let $C$ be the strict transform of $C'$ on $X$.\end{flushleft}

Let $E''$ be a generalized log canonical place of $(X',B'+cC'+M'+cN')$ such that $ga(E'',X',B'+c'C'+M'+c'N')<0$ for any $c'>c$. If $E''\subset X'$, then
\begin{center}
$1=i_1+cj$
\end{center}
for some $i_1\in I_1$ and $j\in J$. Thus by \hyperref[0401]{Lemma 4.1(3)}, $c\in\mathcal{N}_d(I_1,I_2;J)$. Thus, from now on we may assume $E''$ is not contained in $X'$. Let $\psi: X''\rightarrow X'$ be a generalized DLT modification of $(X',B'+cC'+M'+cN')$ such that the center of $E''$ on $X''$ is a divisor. We identify $E''$ with its image on $X''$. Possibly replacing $X$ by resolved data, we may assume $\phi: X\dashrightarrow X''$ is a morphism. Thus, let $B'',C'',M'',N''$ be the pushdown of $B,C,M,N$ from $X$ to $X''$, we have

\begin{center}
$K_{X''}+B''+cC''+M''+cN''=\psi^*(K_{X'}+B'+cC'+M'+cN')$.
\end{center} 

Notice that the coefficients of $C''$ has a one-to-one correspondence with the coefficients of $C'$, and $B''$ is the strict transform of $B'$ plus the reduced exceptional divisor of $\psi$. Thus, $E''$ is a component of $B''$ and $B''\in D(I_1)$. 

Since $ga(E'',X',B'+c'C'+M'+c'N')<0$ for any $c'>c$, we show that there exists a component $F''$ of the exceptional divisor of $\psi$, such that either $C''|_{F''}\not=0$, or $N''|_{F''}\not\equiv_{num}0$. We may write
\begin{center}
$\psi^*N'=N''+E''_1$
\end{center}
and 
\begin{center}
$\psi^*C'=C''+E_2''$.
\end{center}

Since $N$ is nef over $Z$ hence nef over $X'$, $N''$ is nef over $X'$ hence by the negativity lemma, $E_1''$ and $E_2''$ are both effective. Moreover, since $ga(E'',X',B'+c'C'+M'+c'N')<0$ for any $c'>c$, $E''$ has positive coefficient either in $E_1''$ or $E_2''$, and in particular either $E_1''$ or $E_2''$ is not numerically trivial over $X'$. If we cannot find a component $F''$ of $Ex(\psi)$, such that either $C''|_{F''}\not=0$, or $N''|_{F''}\not\equiv_{num}0$, then for any component $F''$ of $Ex(\psi)$, $E''_1|_{F''}\equiv_{num}0$ and $E''_2|_{F''}\equiv_{num}0$. This implies that $E''_1$ and $E_2''$ are numerically trivial over $X'$, which is a contradiction.

Let $\Delta''=B''+cC''$ and $\Gamma''=M''+cN''$, $\Gamma=M+cN$. Let 
\begin{center}
$K_{E''}+\Delta_{E''}+\Gamma_{E''}=(K_{X''}+\Delta''+\Gamma'')|_{E''}$
\end{center}
be the generalized adjunction, then by \hyperref[0306]{Theorem 3.6}, $\Gamma_{E''}=M_{E''}+cN_{E''}$ where $M_{E''}\in\mathcal{NEF}_{E}(I_2)$ and $N_{E''}\in\mathcal{NEF}_{E}(J)$, $\Delta_{E''}\in D(I_1\cup I_2\cup J)=D(I_1\cup J)$. By \hyperref[0307]{Theorem 3.7}, either $N_{E}\not\equiv_{num}0$, or there exists a component of $\Delta_{E''}$ with coefficient of the form
\begin{center}
$\frac{m-1+i+\sum_{p=1}^{l} j_pz_p}{m}$
\end{center}
where $i\in (I_1)_+=I_1$, $l\geq 1$, $j_p\in J\backslash\{0\}$ for each $p$, $z_p\in\mathbb N^+$ for each $p$, $m\in\mathbb N^+$.

Thus $\Delta_{E''}=B_{E''}+C_{E''}$, where $B_{E''}\in D(I_1), C_{E''}\in\mathcal{D}_c(I_1;J)$. Moreover, since $K_{E''}+\Delta_{E''}+\Gamma_{E''}\equiv_{num}0$, we have $c\in\mathcal{N}_d(I_1,I_2;J)$. Thus, $GLCT_{d+1}(I_1,I_2;J)\subset\mathcal{N}_d(I_1,I_2;J)$

Now we show that $GLCT_{d+1}(I_1,I_2;J)\supset\mathcal{N}_d(I_1,I_2;J)$. Thus, by \hyperref[0402]{Lemma 4.2}, we only need to show that $GLCT_{d+1}\supset\mathcal{K}_d(I_1,I_2;J)$. By induction on $d$, we may assume there exists a generalized KLT pair $(X',\Delta'+\Gamma')\in\mathcal{S}_d(I_1,I_2;J,c)$ of dimension $d$. From the definition of $\mathcal{R}_d(I_1,I_2;J,c)$, either $\Delta'\not=0$ or $\Gamma'\not\equiv_{num}0$, thus $K_{X'}$ is not numerical trivial, hence $K_{X'}$ is anti-ample. Since $J$ is closed under addition, all components of $\Delta'$ are of the form
\begin{center}
$\frac{m-1+i+jc}{m}$
\end{center}
for some $i\in I_1$ and $j\in J$. Moreover, we may write $\Gamma=M+cN$ such that $M\in\mathcal{NEF}_{X}(I)$ and $N\in\mathcal{NEF}_X(J)$, and let $M',N'$ be the pushdown of $M,N$ on $X'$, and $\Delta'=B'+C'$ where $B'\in D(I_1)$ and $C'\in\mathcal{D}_c(I_1;J)$.

Now we embed $X'$ into $\mathbb P^r$ for some $r>0$, and consider the affine cone of $\bar X'=C(X',L')\subset\mathbb A^{r+1}$ for the very ample line bundle $L':=\mathscr{O}_{\mathbb P^r}(1)|_{X'}$ of  $X'$. 

Write $B'=\sum b_jB_j', C'=\sum c_jC_j', M'=\sum r_jR_j'$ and $N'=\sum s_jS_j'$ in terms of irreducible components, we let $\bar B_j'$, $\bar C_j'$, $\bar R_j'$ and $\bar S_j'$ be the corresponding  divisors on the cone $\bar X'$, and let $\bar B'=\sum b_j\bar B_j'$,  $\bar C'=\sum c_j\bar C_j'$,  $\bar M'=\sum r_j\bar R_j'$,  and $\bar N'=\sum s_j\bar S_j'$. 

Let $\pi': W'\rightarrow \bar X'$ be the blow-up of the vertex with exceptional divisor $Y'\cong X'$. Let $B_{W'}$, $C_{W'}$, $M_{W'}$ and $N_{W'}$ be the strict transform of $\bar B'$, $\bar C'$, $\bar M'$ and $\bar N'$ on $W'$. Then we have
\begin{center}
$K_{W'}+B_{W'}+C_{W'}+Y'+M_{W'}+cN_{W'}=(\pi')^*(K_{\bar X'}+\bar B'+\bar C'+\bar M'+c\bar N')$.
\end{center}

Moreover, there exists a log resolution $W:=W'\times_{X'}X$. Let $h: W\rightarrow W'$ be the corresponding morphism, and let $M_{W}=h^*M_{W'}$, $N_{W}=h^*N_{W'}$. Then 
$(\bar X',\bar B'+\bar C'+\bar M'+c\bar N')$ is a generalized pair with data $(M_W+cN_W,W\rightarrow \bar X'\rightarrow\mathbb A^1)$. In particular,  $(\bar X',\bar B'+\bar C'+\bar M'+c\bar N')$ is generalized LC but not generalized KLT.

For each component $\bar S'$ of $\bar\Delta'$, if the corresponding coefficient in $\bar B'+\bar C'$ is of the form
\begin{center}
$\frac{m-1+i+jc}{m}$,
\end{center}
for some $i\in I_1$, $j\in J$ and $m=m(S')\in\mathbb N^+$, there is a ramified cover of order $m$ for $\bar S'$. 

In particular, there exists a cover $\zeta: \tilde X'\rightarrow\bar X'$ that ramifies each $\bar S'$ to order $m(S')$. Let $\tilde M'$ and $\tilde N'$ be the pullback of $\bar M'$ and $\bar N'$ to $\tilde X'$. Then we have
\begin{center}
$K_{\tilde X'}+\tilde\Delta'+\tilde M'+c\tilde N'=\zeta^*(K_{\bar X'}+\bar\Delta'+\bar M'+c\bar N')$.
\end{center}
where each component of $\tilde\Delta'$ has the form $i+jc$, where $i\in I_1$ and $j\in J$. Moreover, Let $\tilde W$ be the fiber product of $W$ and $\tilde X'$, and let $\tilde M$, $\tilde N$ be the pullback of $M$ and $N$ to $\tilde W$. Then $\tilde M\in\mathcal{NEF}_{\tilde W}(I_2)$, $\tilde N\in\mathcal{NEF}_{\tilde W}(J)$, and there exists a natural morphism $\tilde W\rightarrow\tilde X'$ such that $\tilde M'$, $\tilde N'$ are the pushdown of $\tilde M$ and $\tilde N$. In particular, this gives $(\tilde X',\tilde\Delta'+\tilde M'+c\tilde N')$ a generalized pair structure with nef part $\tilde M'+c\tilde N'$. Moreover, the generalized log canonical and non-generalized KLT property of $(\tilde X',\tilde\Delta'+\tilde M'+c\tilde N')$ is also preserved.

Since all the coefficients of $\tilde\Delta'$ are of the form $i+jc$, we may write $\Delta'=\tilde B'+c\tilde C'$, such that $\tilde B'\in I_1$, $\tilde C'\in J$.

Now $c=glct(\tilde X',\tilde B'+\tilde M';\tilde C'+\tilde N')\in GLCT_{d}(I_1,I_2;J)$.\hfill$\Box$

\section{The Main Claim}

In this section we prove the main claim of our paper, which is strongly related to \hyperref[0106]{Theorem 1.6}. Instead of working with $GLCT_{d}(I_1,I_2;J)$, we intend to work with $\mathcal{N}_d(I_1,I_2,J)$.

\begin{claim}\label{0501}\rm
Let $d\geq 2$ be an integer, $I_1,I_2,J\subset [0,\infty)$ be three DCC sets, such that $1\in I_1=(I_1)_+$, $I_2\cap [0,1]\subset I_1$, $J$ is closed under addition. Then there exists two DCC  sets $I_1',J'$ contained in $[0,\infty)$ that only depends on $I_1,I_2,J,d$, such that 

(i) $I_1'=(I_1')_+=\bar I_1'$;

(ii) $J'$ is closed under addition;

(iii) $I_1\subset I_1'$, $J\subset J'$;

(v)  All the accumulation points of $\mathcal{N}_d(I_1,I_2;J)$ are contained in either $GLCT_{1}(I_1',\{0\};J')$, or $\mathcal{N}_{d-1}(\bar I_1',I_2,\bar J)$.

Moreover, there exists an integer $N>0$ that only depends on $I_1,I_2,J$ and $d$, such that 

(1) $J'=\frac{1}{N}\bar J$;

(2) $I_1'=(\frac{1}{N}(\bar I_1\cup \bar I_2))_+$.
\end{claim}

For inductive reasons, before we start proving \hyperref[0501]{Claim 5.1}, we prove the following lemma:

\begin{lem}\label{0502}\rm
\hyperref[0106]{Theorem 1.6} holds when $d=2$.
\end{lem}

\begin{flushleft}\textit{Proof}. By \hyperref[0403]{Lemma 4.3}, when $GLCT_{2}(I_1,I_2;J)=\mathcal{N}_1(I_1,I_2;J)$. For any $c\in\mathcal{N}_1(I_1,I_2;J)$, $c$ is a solution of the equation
\begin{center}
$2=\sum_k cj_k+\sum_k i_{1,k}+\sum_k i_{2,k}$
\end{center}
where each $j_k\in J$, $i_{1,k}\in I_1$ and $i_{2,k}\in I_2$. Thus, the only accumulation point $\tilde c$ of all the possible $c\in\mathcal{N}_1(I_1,I_2;J)$ is either $0\in GLCT_{1}(I_1,0,J)$, or a solution of the equation
\begin{center}
$2=\sum_k \tilde c\tilde j_k+\sum_k \tilde i_{1,k}+\sum_k \tilde i_{2,k}$
\end{center}
where each $\tilde j_k\in \bar J$, $\tilde i_{1,k}\in \bar I_1$ and $\tilde i_{2,k}\in \bar I_2$.\end{flushleft}

Let $J'=\frac{1}{2}\bar J$, $I_1'=(\frac{1}{2}(\bar I_1\cup\bar I_2))_+$. Notice that since $J$ is closed under addition, so is $\bar J$ and $J'$, thus $u:=\sum\frac{\tilde j_k}{2}\in J'$ and $v:=\frac{1}{2}(\sum_k \tilde i_{1,k}+\sum_k \tilde i_{2,k})\in I_1'$. Now let $X'=X=\mathbb P^1$, $p\in X$ be a point, $B=vp$ and $C=up$, $M=N=0$. Then
\begin{center}
$c=glct(X',B;C)$.
\end{center}

Thus $c\in GLCT_1(I_1',0;J')$.\hfill$\Box$

\begin{flushleft}\textit{Proof of \hyperref[0501]{Claim 5.1}}.\end{flushleft}

\begin{flushleft}\textit{Step 1}. (Basic Setting)\end{flushleft}

Fix $d$ and pick $I_1,I_2,J$ as in the assumptions of \hyperref[0501]{Claim 5.1}. Since $\mathcal{N}_d(I_1,I_2;J)=GLCT_{d+1}(I_1,I_2;J)$ under our assumptions, and by \hyperref[0303]{Theorem 3.3}, $GLCT_{d+1}(I_1,I_2;J)$ is an ACC set, $\mathcal{N}_d(I_1,I_2;J)$ is an ACC set. Thus for any accumulation point $c$ of $\mathcal{N}_d(I_1,I_2;J)$, we may assume $c$ is given by a strict decreasing sequence of numbers $c_1,\dots,c_i,\dots\in\mathcal{N}_d(I_1,I_2;J)$ that converges to $c$. From the definition of $\mathcal{N}_d(I_1,I_2;J)$, for each $i$, there exists a projective generalized pair $(X_i',\Delta_i'+\Gamma_i')$ with data $(\Gamma_i,X_i\xrightarrow{f_i}X_i')$ of dimension $d$, such that:

(i) $K_{X_i'}+\Delta_i'+\Gamma_i'\equiv_{num}0$;

(ii) $\Gamma_i=M_i+c_iN_i$, where $M_i\in\mathcal{NEF}_{X_i}(I_2)$, $N_i\in\mathcal{NEF}_{X_i}(J)$;

(iii) $M_i':=(f_i)_*M_i$, $N_i':=(f_i)_*N_i$;

(iv) $\Delta_i'=B_i'+C_i'$, such that $B_i'\in D(I_1)$, and $C_i'\in\mathcal{D}_{c_i}(I_1;J)$.

(v) Either $C_i'\not=0$, or $N_i\not\equiv_{num}0$.

By \hyperref[0402]{Lemma 4.2} and using induction on dimensions, then using \hyperref[0401]{Lemma 4.1}, we may also assume the following:

(vi) $X_i'$ is $\mathbb Q$-factorial, $\rho(X_i')=1$, $(X_i',\Delta_i'+\Gamma_i')$ is generalized KLT.

Since $0\in\mathcal{N}_{d-1}(I_1,I_2;J)$ is clear, we may assume $c>0$. 

Let $a_i=gmld(X_i',\Delta_i'+\Gamma'_i)$. We may assume $a_i$ either forms a decreasing or an increasing sequence (although if we assume ACC for generalized minimal log discrepancies, $a_i$ must be decreasing, but this is not important for our proof). Let $a=lim\ a_i$. 

For every real number $0<\epsilon<1$, we may rewrite $B_i'=B_{i,\epsilon+}'+B_{i,\epsilon-}'$ and $C_i'=C_{i,\epsilon+}'+C_{i,\epsilon-}'$, such that 

1) $B_{i,\epsilon+}'\wedge B_{i,\epsilon-}'=0$ and $C_{i,\epsilon+}'\wedge C_{i,\epsilon-}'=0$;

2) All the coefficients of $B_{i,\epsilon-}'$ and $C_{i,\epsilon-}'$ are $\geq 1-\epsilon$, and all the coefficients of $B_{i,\epsilon+}'$ and $C_{i,\epsilon+}'$ are $<1-\epsilon$.

We define $A_{i,\epsilon}'=B_{i,\epsilon-}'+C_{i,\epsilon-}'$.

For every component $V_i'$ of $C_i'$, the coefficient of $V_i'$ in $C_i'$ is of the form
\begin{center}
$s_{V'_i}=\frac{m-1+i_1+jc_i}{m}$
\end{center}
where $i_1\in I_1$, $j\in J$, $m\in\mathbb N^+$. For each $s_{V'_i}$, any real number $0\leq t\leq c_i$, we define $s(t)_{V'_i}=\frac{m-1+i_1+jt}{m}$. For any $0<\epsilon<1$, we define $C'(t)_{i,\epsilon+}=\sum_{V'_i\ is\ a\ component\ of\ C_{i,\epsilon+}'}s(t)_{V'_i}V_i'$.

\begin{flushleft}\textit{Step 2}. We prove the following two lemmas, and construct a sufficiently small $\epsilon$ with good properties.\end{flushleft}

\begin{lem}\label{0503}\rm
Assumptions as in \textit{Step 1}. $c$ is contained in a DCC set that only depends on $I_1,I_2$ and $J$. In particular, for any component $V_i'$ of $C_i'$, $s(c)_{V'_i}$ is contained in a DCC set that only depends on $I_1,I_2,J$. 
\end{lem}

\begin{flushleft}\textit{Proof}.  First we show that $c$ is contained in a DCC set that only depends on $I_1,I_2$ and $J$. If not, suppose that there exists a strict increasing sequence $\{c^j\}_{j=1}^{\infty}$ where each $c^j$ is the accumulation point of $\{c_i^j\}_{i=1}^{\infty}$ from above, and where $c_i^j\in\mathcal{N}_d(I_1,I_2;J)$ for every $i,j$.\end{flushleft}

Now for every $k\geq 1$, we inductively pick $i_k$ such that $c^k_{i_k}<c^{k+1}$. In this way
\begin{center}
$\{c^k_{i_k}\}_{k=1}^{\infty}$
\end{center}
is a strictly increasing sequence, which is not possible.

Notice that $s(c)_{V'_i}\in D(I_1\cup c\cdot J)$ which is a DCC set. Since we have already shown that $c$ is contained in a DCC set that  only depends on $I_1,I_2$ and $J$, $D(I_1\cup c\cdot J)$ only depends on $I_1,I_2$ and $J$.\hfill$\Box$

\begin{lem}\label{0504}\rm
For any DCC set $I, J\subset [0,\infty)$, and integer $d>0$, there exists $1>\epsilon=\epsilon(I,J)>0$ that only depends on $I,J,d$ that satisfies the following properties: 

For any generalized log canonical pair $(X',A'+B'+M')$ of dimension $\leq d$ with data $(M,X\rightarrow X')$, such that $M'\in\mathcal{NEF}_{X}(J)$, $B'\in I$, and $A'\in [1-\epsilon, 1]$, let $\bar A'=$Supp$A'$. Then:

(1) $(X',\bar A'+B'+M')$ is generalized log canonical.

(2) If in addition, $X'$ is of Picard number $1$, $-(K_{X'}+A'+B'+M')$ is pseudoeffective and $K_{X'}+\bar A'+B'+M'$ is pseudoeffective, then 
\begin{center}
$K_{X'}+\bar A'+B'+M'\equiv_{num}0$.
\end{center}

(3) If $A'=(1-\alpha)\bar A'$ for some $\alpha\geq 0$, and $K_{X'}+A'+B'+M'$ is numerically trivial, then $\alpha=0$.

\end{lem}

\begin{flushleft}\textit{Proof}. We may only care about generalized pairs $(X',A'+B'+M')$ of dimension $d$.\end{flushleft}

To prove (1), since $GLCT_{d}(I,J;\{0,1\})$ is an ACC set by \hyperref[0303]{Theorem 3.3}, there exists 
\begin{center}
$r=\underset{t<1}{sup}\{t\in GLCT_{d}(I,J;\{0,1\})\}<1$.
\end{center}

We may pick $\epsilon<\frac{1-r}{2}$. Since $(1-\epsilon)\bar A'\leq A'$, $(X',(1-\epsilon)\bar A'+B'+M')$ is generalized log canonical. Thus,
\begin{center}
$1\geq glct(X',B'+M';\bar A')\geq 1-\epsilon>r$.
\end{center}

Thus $glct(X',B'+M';\bar A')=1$, and hence  $(X',\bar A'+B'+M')$ is generalized log canonical.

To prove (2), since $X'$ is of Picard number $1$, there exists $\delta\geq 0$ such that
\begin{center}
$K_{X'}+(1-\delta)\bar A'+B'+M'\equiv_{num}0$
\end{center}

Since $-(K_{X'}+A'+B'+M')$ is pseudoeffective, $-(K_{X'}+(1-\epsilon)\bar A'+B'+M')$ is pseudoeffective, thus
\begin{center}
$1-\epsilon\leq 1-\delta\leq 1$.
\end{center}

If we cannot find such an $\epsilon$, we may construct a strict decreasing sequence $\epsilon_i$ that converges to $0$, a sequence of generalized pairs in the following way: $(X'_i,A'_i+B'_i+M'_i)$ is a generalized pair, such that $X'_i$ is of Picard number $1$, $-(K_{X_i'}+A_i'+B_i'+M_i')$ is pseudoeffective and $K_{X_i'}+\bar A_i'+B_i'+M_i'$ is pseudoeffective but not numerically trivial, all the coefficients of $A_i'$ are $\geq 1-\epsilon_i$. Let $\delta_i$ be the unique number such that 
\begin{center}
$K_{X_i'}+(1-\delta_i)\bar A_i'+B_i'+M_i'\equiv_{num}0$
\end{center}

Then $1-\epsilon_i<1-\delta_i<1$.

Since $1-\epsilon_i$ converges to $1$, $1-\delta_i$ converges to $1$. Thus $\{1-\delta_i\}$ is not ACC, which contradicts to \hyperref[0304]{Theorem 3.4}.

If (3) does not hold, then there exists a strictly decreasing sequence of $\epsilon_i$ that converges to $0$, a sequence of numbers $\alpha_i>0$, a sequence of generalized pairs in the following way: $(X'_i,A'_i+B'_i+M'_i)$ is a generalized pair, such that $A'=(1-\alpha_i)\bar A'\in [1-\epsilon_i, 1]$, and $K_{X'_i}+A'_i+B'_i+M'_i\equiv_{num}0$. Thus, $1-\epsilon_i\leq 1-\alpha_i<1$, hence after possibly passing to a subsequence, we may assume $1-\alpha_i$ is strictly increasing. But $\{1-\alpha_i\}$ is not ACC, which contradicts to \hyperref[0304]{Theorem 3.4}.\hfill$\Box$

\begin{flushleft}\textit{Step 3}. By \hyperref[0503]{Lemma 5.3}, any numbers of the form 
\begin{center}
$\frac{m-1+i_1+jc}{m}$
\end{center}
where $m\in\mathbb N^+$, $i_1\in I_1$ and $j\in J$ are contained in a DCC set $\mathcal{P}$ that only depends on $I_1$, $I_2$ and $J$, and any numbers of the form
\begin{center}
$i_2+cj$
\end{center}
are contained in another DCC set $\mathcal{Q}$ that only depends on $I_1,I_2$ and $J$. We may assume that $\mathcal{P},\mathcal{Q}$ both contain $I_1, I_2, J$, and in particular they both contain $1$.\end{flushleft}

Thus, we may define $\epsilon=\epsilon(\mathcal{P},\mathcal{Q})$ as in \hyperref[0504]{Lemma 5.4}.

Now we consider $a$ (which is defined in \textit{Step 1}). Then either $a<\epsilon$, or $a\geq\epsilon$. If $a<\epsilon$, after possibly passing to a subsequence we may assume that for each $i$, $a_i<\epsilon$. Thus, there are four cases:

\begin{flushleft}\textbf{Case 1}\label{Case1} $a<\epsilon$ and $A'_{i,\epsilon}\not=0$;

\textbf{Case 2}\label{Case2} $a<\epsilon$ and $A'_{i,\epsilon}=0$, while $a_i>0$ for every $i$;

\textbf{Case 3}\label{Case3} after possibly passing to a subsequence, $a_i=0$ for every $i$;

\textbf{Case 4}\label{Case4} $a\geq\epsilon$.\end{flushleft}

We shall deal with \hyperref[Case1]{\textbf{Case 1}} in \textit{Step 4}, \hyperref[Case2]{\textbf{Case 2}} in \textit{Step 5}, and \hyperref[Case3]{\textbf{Case 3}} in \textit{Step 6} and \hyperref[Case4]{\textbf{Case 4}} in \textit{Step 7}.

\begin{flushleft}\textit{Step 4}. In this step we deal with \hyperref[Case1]{\textbf{Case 1}}. First we prove the following lemma:\end{flushleft}

\begin{lem}\label{0505}\rm
Assumptions as in \textit{Step 4}. We are done with the case when $C'_{i,\epsilon+}=0$ for every $i$ and $N_i'\equiv_{num}0$ for every $i$.
\end{lem}

\begin{flushleft}\textit{Proof}. If not, passing to a subsequence we may assume that $C'_{i,\epsilon+}=0$ for every $i$ and $N_i'\equiv_{num}0$ for every $i$. Thus
\begin{center}
$K_{X_i'}+B'_i+C'_{i,\epsilon-}+M'_i\equiv_{num}0$.
\end{center}\end{flushleft}

Notice that all the coefficients of $C'_{i,\epsilon-}$ are $\geq 1-\epsilon$, $B_i'\in\mathcal{P}$ and $M'_i\in\mathcal{NEF}_{X_i}(\mathcal{Q})$. Thus, \hyperref[0504]{Lemma 5.4} tells us that 
\begin{center}
$K_{X_i'}+B'_i+\bar C'_{i,\epsilon-}+M'_i\equiv_{num}0$
\end{center}
where $\bar C'_{i,\epsilon-}=$Supp$C'_{i,\epsilon-}$.

This tells us that $\bar C'_{i,\epsilon-}=C'_{i,\epsilon-}$. Since $N'_i\equiv_{num}0$, $C_i'\not\equiv_{num}0$, and so $C'_{i,\epsilon-}=C_i'\not=0$ in this case. Hence, from the construction of $C_i'$, there exists $u_{1,i}\in I_1$ and $j_i\in J$ such that
\begin{center}
$c_i=\frac{1-u_{1,i}}{j_i}$
\end{center}
we may assume that $j_i$ converges to $j$ and $u_{1,i}$ converges to $u_1$. Then $j\in J$ since $J=\bar J$, and $u_1\in \bar I_1$ which is a DCC set that only depends on $I_1$. Thus $c_i\in\mathcal{N}_1(\bar I_1,I_2;J)\subset\mathcal{N}_{d-1}(\bar I_1,I_2;\bar J)$. Thus \hyperref[0501]{Claim 5.1} holds in this case.\hfill$\Box$
 
Thus from now on until the end of \textit{Step 4}, we may assume that either $N_i'\not\equiv_{num}0$ for every $i$, or $C'_{i,\epsilon+}\not=0$ for every $i$.

By our assumptions, we have
\begin{center}
$K_{X_i'}+A'_i+B'_{i,\epsilon+}+C'_{i,\epsilon+}+M'_i+cN'_i\equiv_{num}0$.
\end{center}

Let $\bar A'_i=$Supp$A'_i$. By our construction of $\mathcal{P},\mathcal{Q}$ and $\epsilon$, we have

\begin{center}
$(X_i',\bar A'_i+B'_{i,\epsilon+}+C'_{i,\epsilon+}+M'_i+cN'_i)$
\end{center}
is generalized log canonical. 

Let $S_i'$ be an irreducible component of $A'_i$. Notice that since $X_i'$ is of Picard number $1$, if $N'_i\not\equiv_{num}0$, $N'_i|_{S_i'}\not\equiv_{num}0$, and if $C'_i\not=0$, $C'_i|_{S_i'}\not=0$. 

We consider the following thresholds: we define
\begin{center}
$r_i=\underset{t>0}{sup}\{K_{X_i'}+\bar A'_i+B'_{i,\epsilon+}+C'(t)_{i,\epsilon+}+M'_i+tN'_i$ is pseudoeffective$\}$.
\end{center}

Since $X_i'$ is of Picard number $1$, $r_i$ is the unique number such that 
\begin{center}
$K_{X_i'}+\bar A'_i+B'_{i,\epsilon+}+C'(r_i)_{i,\epsilon+}+M'_i+r_iN'_i\equiv_{num}0$.
\end{center}

By \hyperref[0304]{Theorem 3.4}, since either $C'(r_i)_{i,\epsilon+}\not=0$ or $N'_i\not\equiv_{num}0$,  $r_i$ must be contained in an ACC set. Thus after possibly passing to a subsequence, we may assume that $r_i$ is a decreasing sequence and converges to $r$. Notice that $r_i\leq c_i$, we must have $r\leq c$.

Thus, there are three cases:

\begin{flushleft}\textbf{Case 1.1}\label{Case0101} $r<c$. In this case, after possibly passing to a subsequence, we may assume that $r_i<c$ for every $i$. Thus, 
\begin{center}
$K_{X_i'}+\bar A'_i+B'_{i,\epsilon+}+C'(c)_{i,\epsilon+}+M'_i+cN'_i$
\end{center}
is pseudoeffective but not numerically trivial. Now we define
\begin{center}
$t_i=\underset{t>0}{sup}\{K_{X_i'}+t_i\bar A'_i+B'_{i,\epsilon+}+C'(c)_{i,\epsilon+}+M'_i+cN'_i$ is pseudoeffective$\}$.
\end{center}\end{flushleft}

Since $X_i'$ is of Picard number $1$, $t_i$ is the unique number such that 
\begin{center}
$K_{X_i'}+t_i\bar A'_i+B'_{i,\epsilon+}+C'(c)_{i,\epsilon+}+M'_i+cN'_i\equiv_{num}0$.
\end{center}

By our construction, $t_i<1$. But notice that
\begin{center}
$K_{X_i'}+A'_i+B'_{i,\epsilon+}+C'(c_i)_{i,\epsilon+}+M'_i+c_iN'_i\equiv_{num}0$.
\end{center}
and 
\begin{center}
$A'_i\geq (1-\epsilon)\bar A'_i$
\end{center}

We have
\begin{center}
$t_i>1-\epsilon$.
\end{center}

Thus, by our construction of $\mathcal{P},\mathcal{Q}$ and $\epsilon$, we get a contradiction. This finishes \hyperref[Case0101]{Case 1.1}.

\begin{flushleft}\textbf{Case 1.2}\label{Case0102} After possibly passing to a subsequence, $r=c$ and $r_i\not=c$ for every $i$.\end{flushleft}

In this case, we pick a component $S_i'$ of $\bar A_i'$ and consider the generalized adjunction of 
\begin{center}
$K_{X_i'}+\bar A'_i+B'_{i,\epsilon+}+C'(r_i)_{i,\epsilon+}+M'_i+r_iN'_i$
\end{center}
to $S_i'$.

Since either $C'(r_i)_{i,\epsilon+}\not=0$ or $N'_i\not\equiv_{num}0$, by \hyperref[0307]{Theorem 3.7}, we have $r_i\in\mathcal{N}_{d-1}(I_1,I_2;J)$. Since $c$ is the accumulation point of $\{r_i\}$, the claim is followed by induction on dimensions in this case. This finishes \hyperref[Case0102]{Case 1.2}.

\begin{flushleft}\textbf{Case 1.3}\label{Case0103} After possibly passing to a subsequence, $r_i=c$ for every $i$. In this case, we consider the generalized adjunction of 
\begin{center}
$K_{X_i'}+\bar A'_i+B'_{i,\epsilon+}+C'(c)_{i,\epsilon+}+M'_i+cN'_i$
\end{center}
to $S_i'$, and we find that $c\in\mathcal{N}_{d-1}(I_1,I_2;J)$. Moreover, it is clear that under either assumption of \hyperref[0501]{Lemma 5.1}, the corresponding claim holds. This finishes \hyperref[Case0103]{Case 1.3}, and thus finishes \hyperref[Case1]{Case 1}.\end{flushleft}

Since for all other cases, $C'_{i,\epsilon+}=C'_i$, for simplicity we define $C'(t)_i=C'(t)_{i,\epsilon+}$ for any real number $0\leq t\leq c_i$ in the rest of the proof.

\begin{flushleft}\textit{Step 5}. In this step we deal with \hyperref[Case2]{Case 2}. Since $a_i>0$ for every $i$ and $a<\epsilon$, after possibly passing to a subsequence, we may assume $a_i<\epsilon$ for every $i$, and we let $\pi_i: X_i''\rightarrow X_i'$ be the extraction of an irreducible reduced divisor $E_i''$ such that
\begin{center}
$ga(E_i'',X_i',B'_i+C'_i+M'_i+c_iN'_i)=a_i$.
\end{center}\end{flushleft}

After possibly replacing $X_i$, we may assume that $X_i\dashrightarrow X_i''$ is a morphism. We let $B_i'',C_i'',C''(t)_i$ be the strict transforms of $B_i', C_i', C'(t)_i$ on $X_i''$, and let $M''_i$, $N''_i$ be the pushdown of $M_i$, $N_i$ to $X_i''$. 

We have
\begin{center}
$K_{X_i''}+(1-a_i)E''_i+B''_i+C''_i+M''_i+c_iN_i''=\pi_i^*(K_{X_i'}+B'_i+C'_i+M'_i+c_iN'_i)$.
\end{center}

Since $K_{X_i'}+B'_i+C'_i+M'_i+c_iN'_i\equiv_{num}0$, $K_{X_i''}+B''_i+C''_i+M''_i+c_iN_i''$ is not pseudoeffective, and we may run a $(K_{X_i''}+B''_i+C''_i+M''_i+c_iN_i'')$-MMP with scaling. Let $\phi_i$ be any partial MMP of this minimal model program. Then $\phi_i$ is either birational, or a Mori fiber space.

If $\phi_i$ is a birational map for every $i$, suppose $\phi_i: X_i''\dashrightarrow X_i'''$ is the map. We let $B'''_i, C'''_i, C'''(t)_i, E'''_i$ be the birational transform of $B''_i,C''_i, C''(t)_i$ and $E''_i$ on $X'''_i$. Possibly replacing $X$ by a resolved data, we may assume that $X_i\dashrightarrow X_i'''$ is a morphism, and we define $M_i'''$ and $N_i'''$ to be the pushdown of $M_i$ and $N_i$ to $X'''_i$.

Notice since $\phi_i$ is $E_i''$-positive, $E_i''$ is not contracted. If $C''_i$ is contracted and $N_i'''\equiv_{num}0$, we have
\begin{center}
$K_{X'''_i}+(1-a_i)E_i'''+B'''_i+M_i'''\equiv_{num}0$.
\end{center}

But this contradicts to \hyperref[0504]{Lemma 5.4(3)}. Thus, either $C''_i$ is not contracted, or $N_i'''\not\equiv_{num}0$. Since $\phi_i$ can be any partial MMP, and the $(K_{X_i''}+B''_i+C''_i+M''_i+c_iN_i'')$-MMP with scaling always terminates with a Mori fiber space, we may assume that there exists a Mori fiber space structure $\psi_i: X_i'''\rightarrow T_i$. 

If $\phi_i$ is not a birational map except finitely many $i$, we may assume $\phi_i:X_i''\rightarrow T_i$ is not a birational map for every $i$, and we define $\psi_i=\phi_i$. For simplicity, we may write $X_i'''=X_i'',B'''_i=B''_i, C'''_i=C''_i, C'''(c)_i=C''(c)_i, E'''_i=E''_i, M'''_i=M''_i$ and $N'''_i=N''_i$ in the rest of the \textit{Step 4}.

\begin{flushleft}\textit{Step 5.1}. If $dim T_i>0'''$, let $F'''_i$ be the general fiber of $\psi_i$ and we consider the projective generalized pair
\begin{center}
$(F_i''',(1-a_i)E'''_i|_{F_i'''}+B'''_i|_{F_i'''}+C'''_i|_{F_i'''}+M'''_i|_{F_i'''}+c_iN_i'''|_{F_i'''})$
\end{center}\end{flushleft}

Now we are in the lower dimension case of \textit{Step 5} after replacing $X_i'',E_i'',B''_i, C''_i,M''_i$ and $N''_i$ by $F_i''',E'''_i|_{F_i'''}, B'''_i|_{F_i'''}, C'''_i|_{F_i'''}, M'''_i|_{F_i'''}$ and $N_i'''|_{F_i'''}$. The claim is followed by induction on dimensions in this case. 

\begin{flushleft}\textit{Step 5.2}. If $T_i$ is a point, then $X_i'''$ is of Picard number $1$. By \hyperref[0504]{Lemma 5.4(1)}, we have
\begin{center}
$(X'''_i,E_i'''+B'''_i+C'''(c)_i+M'''_i+cN'''_i)$
\end{center}
is log canonical. Now we define
\begin{center}
$r_i=\underset{t>0}{inf}\{K_{X'''_i}+E_i'''+B'''_i+C'''(t)_i+M_i'''+tN_i'''$ is pseudoeffective$\}$
\end{center}\end{flushleft}

Since $X'''_i$ is of Picard number $1$, it is clear that $r_i$ is the unique number such that 
\begin{center}
$K_{X'''_i}+E_i'''+B'''_i+C'''(r_i)_i+M_i'''+r_iN_i'''\equiv_{num}0$.
\end{center}

In particular, $r_i<c_i$ since either $C'''(r_i)_i\not=0$ or $N_i'''\not\equiv_{num}0$. Moreover, if $r_i$ is not contained in an ACC set, it contradicts to \hyperref[0304]{Theorem 3.4}. Thus, after possibly passing to a subsequence, we may assume that $r_i$ is decreasing, and $r=lim\ r_i$.

\begin{flushleft}\textbf{Case 2.1}\label{Case0201}. $r<c$. In this case, after possibly passing to a subsequence, we may assume that $r_i<c$ for every $i$. Thus, 
\begin{center}
$K_{X_i'''}+E'''_i+B''_i+C'''(c)_{i}+M'''_i+cN'''_i$
\end{center}
is pseudoeffective but not numerically trivial. Now we define
\begin{center}
$t_i=\underset{t>0}{sup}\{K_{X_i'''}+t_iE'''_i+B'''_{i}+C'''(c)_{i}+M'''_i+cN'''_i$ is pseudoeffective$\}$.
\end{center}\end{flushleft}

Since $X_i'''$ is of Picard number $1$, $t_i$ is the unique number such that 
\begin{center}
$K_{X_i'''}+t_iE'''_i+B'''_{i}+C'''(c)_{i}+M'''_i+cN'''_i\equiv_{num}0$.
\end{center}

By our construction, $t_i<1$. But notice that
\begin{center}
$K_{X_i'''}+(1-a_i)E'''_i+B'''_i+C'''_i+M'''_i+c_iN_i'''\equiv_{num}0$.
\end{center}

We have
\begin{center}
$t_i>1-a_i$.
\end{center}

Thus, by our construction of $\mathcal{P},\mathcal{Q}$ and $\epsilon$, we get a contradiction. This finishes \hyperref[Case0201]{Case 2.1}.

\begin{flushleft}\textbf{Case 2.2}\label{Case0202}. After possibly passing to a subsequence, $r=c$ and $r_i\not=c$ for every $i$.\end{flushleft}

In this case, we consider the generalized adjunction of 
\begin{center}
$K_{X_i'''}+E'''_i+B'''_{i}+C'''(r_i)_{i}+M'''_i+r_iN'''_i$
\end{center}
to $E_i'''$.

Since either $C'''(r_i)_{i}\not=0$ or $N'''_i\not\equiv_{num}0$, by \hyperref[0307]{Theorem 3.7}, we have $r_i\in\mathcal{N}_{d-1}(I_1,I_2;J)$. Since $c$ is the accumulation point of $\{r_i\}$, the claim is followed by induction on dimensions in this case. This finishes \hyperref[Case0202]{Case 2.2}.

\begin{flushleft}\textbf{Case 2.3}\label{Case0203}. After possibly passing to a subsequence, $r_i=c$ for every $i$. In this case, we consider the generalized adjunction of 
\begin{center}
$K_{X_i'''}+E'''_i+B'''_{i}+C'''(c)_{i}+M'''_i+cN'''_i$
\end{center}
to $E_i'''$, and we find that $c\in\mathcal{N}_{d-1}(I_1,I_2;J)$. Moreover, it is clear that under either assumption of \hyperref[0501]{Lemma 5.1}, the corresponding claim holds. This finishes \hyperref[Case0203]{Case 2.3}, and thus finishes \hyperref[Case2]{Case 2}.\end{flushleft}

\begin{flushleft}\textit{Step 6}. In this case we deal with \hyperref[Case3]{Case 3}. After possibly passing to a subsequence, we may assume that $a_i=0$ for every $i$. Thus, $(X_i',\Delta'+\Gamma')$ is generalized log canonical but not generalized KLT, and we may let $\pi_i: X_i''\rightarrow X_i'$ be a generalized DLT modification of $(X_i',\Delta'+\Gamma')$.\end{flushleft} 

Possibly replacing $X_i$ by a resolved data, we may assume that $X_i\rightarrow X_i''$ is s morphism. Let $B''_i, C''_i, C''(t)_i$ be the strict transform of $B'_i, C_i$ and $C'(t)_i$ on $X_i''$, $M_i''$ and $N_i''$ be the pushdown of $M_i$ and $N_i$ to $X_i''$, and let $E_i$ be the reduced exceptional divisor of $\pi_i$. Then we have
\begin{center}
$K_{X_i''}+E_i''+B''_i+C''_i+M''_i+c_iN''_i=\pi_i^*(K_{X_i'}+B'_i+C'_i+M'_i+c_iN')$.
\end{center}

Since $K_{X_i''}+E_i''+B''_i+C''_i+M''_i+c_iN''_i\equiv_{num}0$, $K_{X_i''}+B''_i+C''_i+M''_i+c_iN''_i$ is not pseudoeffective, and we may run a $(K_{X_i''}+B''_i+C''_i+M''_i+c_iN''_i)$-MMP with scaling of some ample divisor. Let $\phi_i$ be any partial MMP. After possibly replacing $X_i$ by a resolved data, we may assume $X_i\rightarrow X_i'''$ is a morphism.

If $\phi_i$ is birational,  we suppose $\phi_i: X_i''\dashrightarrow X_i'''$ is the map. We let $E_i''',B_i''',C'''_i$ be the strict transform of $E_i'',B_i'',C''_i$ on $X_i'''$, and $M_i''',N_i'''$ be the pushdown of $M_i,N_i$ to $X_i'''$. Notice that $\phi_i$ is $E_i''$-positive, so $E_i''$ is not contracted.

If $C_i'''=0$ and $N_i'''\equiv_{num}0$ for every $i$ and $C''_i\not=0$ for every $i$, then $C''_i$ is contracted. In this case, $C_i''$ intersects $E_i''$, and in particular we may pick a component $S_i''$ of $E_i''$ such that $C_i''|_{E_i''}\not=0$. Now we consider the generalized adjunction of $K_{X_i''}+E_i''+B''_i+C''_i+M''_i+c_iN''_i$ to $S_i''$. By \hyperref[0307]{Theorem 3.7}, we have $c_i\in\mathcal{N}_{d-1}(I_1,I_2;J)$ and we are done.

If $C_i'''=0$, $N_i'''\equiv_{num}0$ for every $i$, and after possibly passing to a subsequence, $C''_i=0$ for every $i$, then $N_i''\not\equiv_{num}0$.

Let $p_i: X_i\rightarrow X_i''$ and $q_i: X_i\rightarrow X_i'''$ be the two morphisms. Since $\phi_i$ does not extract any divisors,  we can write
\begin{center}
$p_i^*N_i''=q_i^*N_i'''+F_i$.
\end{center}
where $F_i$ is an $\mathbb R$-divisor that is exceptional over $X_i'''$. Since $N_i'''\equiv_{num}0$, $p_i^*N_i\equiv_{num}F_i$ is exceptional over $X_i''$, hence by the negativity lemma, $-F_i\geq 0$. However, since $(p_i)_*F_i=0$, $\phi_i$ is $N_i''$-non-negative, hence $F_i\geq 0$. Thus, $F_i=0$. Pick a general curve $\Sigma''_i$ on $X_i''$ such that $N_i''\cdot\Sigma_{i}''>0$, and let $\Sigma_i$ be its strict transform on $X_i$. Then we have
\begin{center}
$0<N_i''\cdot\Sigma_i''=p_i^*N_i''\cdot\Sigma_i=q_i^*N_i'''\cdot\Sigma_i=N_i'''\cdot (q_i)_*\Sigma_i=0$
\end{center}
which is not possible.

Thus, for the rest of \textit{Step 6} we may assume that either $C_i'''\not=0$ or $N_i'''\equiv_{num}0$ for every $i$. Since $\phi_i$ can be any partial MMP, and since in this case all MMP with scaling terminates with a Mori fiber space, we may assume there exists a Mori fiber space structure $\psi_i: X'''_i\rightarrow T_i$. 

If $\phi_i$ itself is a Mori fiber space, for simplicity below, in the rest of \textit{Step 6} we may write $X'''_i=X_i,B'''_i=B_i'', C_i'''=C_i'', E_i'''=E_i'', M_i'''=M_i''$ and $N_i'''=N_i''$, and define $\psi_i=\phi_i: X'''_i\rightarrow T_i$.

If $dim T_i>0$, let $F'''_i$ be the general fiber of $\psi_i$ and we consider the projective generalized pair
\begin{center}
$(F_i''',E'''_i|_{F_i'''}+B'''_i|_{F_i'''}+C'''_i|_{F_i'''}+M'''_i|_{F_i'''}+c_iN_i'''|_{F_i'''})$
\end{center}

We have $(F_i''',E'''_i|_{F_i'''}+B'''_i|_{F_i'''}+C'''_i|_{F_i'''}+M'''_i|_{F_i'''}+c_iN_i'''|_{F_i'''})\in\mathcal{R}_{d-dim T_i}(I_1,I_2;J,c_i)$, thus $c_i\in\mathcal{N}_{d-dim T_i}(I_1,I_2;J)$. The claim is followed by induction on dimensions in this case. 

If $dim T_i=0$, $X_i'''$ is of Picard number $1$. The theorem is followed by \hyperref[Case1]{Case 1}. 

This finishes \hyperref[Case3]{Case 3}.

\begin{flushleft}\textit{Step 7}. In this case we deal with \hyperref[Case4]{Case 4}, and hence finish the proof.\end{flushleft}

In this case, $a\geq\epsilon$, so after possibly passing to a subsequence, we may assume that $a_i\geq\frac{\epsilon}{2}$ for every $i$. Now for every $i$, $(X_i',\Delta_i'+\Gamma_i')$ is $\frac{\epsilon}{2}$-generalized log canonical and hence in particular, $X_i'$ is an $\frac{\epsilon}{2}$-log canonical Fano variety. Thus, $X_i'$ is contained in a bounded family that only depends on $\epsilon$, which indeed, only depends on $I_1,I_2,J$ and $d$ by \hyperref[0311]{Theorem 3.11}.

We define $R_0= \lceil\frac{2}{\epsilon}\rceil!$. Notice that for any component $V_i'$ of $B_i'$ or $C_i'$, if the corresponding coefficient equals to 
\begin{center}
$\frac{w-1+u+vc_i}{w}$
\end{center}
where $w\in\mathbb N^+$, $v\in J$ and $u\in I_1$, then $w\leq\lceil\frac{2}{\epsilon}\rceil$.

Possibly passing to a subsequence, we may assume that the number of nonzero coefficients of $B'_i$, $C'_i$ are fixed, and let them be $p_b$, $p_c$; if we write $M_i=\sum m_{i,j}M_{i,j}$ and $N_i=\sum n_{i,j}N_{i,j}$, where each for each $i,j$, $m_{i,j}\not=0$ implies $M_{i,j}\not\equiv_{num}=0$ and $n_{i,j}\not=0$ implies $N_{i,j}\not\equiv_{num}0$, we may assume that the number of $m_{i,j}$ for each $i$ and the number of $n_{i,j}$ for each $i$ are also both fixed, and let them be $p_m$ and $p_n$.

We write $B'_i=\sum_{j=1}^{p_b} b_{i,j}B'_{i,j}$ and $C'_i=\sum_{j=1}^{p_c}c_{i,j}C'_{i,j}$ in terms of their irreducible components. Possibly passing to a subsequence again, we may assume that
\begin{center}
$b_{i,j}=\frac{w_j-1+u_{i,j}}{w_j}$
\end{center}
where $u_{i,j}\in I_1$, and $w_j\in\mathbb N^+$ does not depend on $i$. We may assume 

\begin{center}
$c_{i,j}=\frac{w'_j-1+u'_{i,j}+c_iv_{i,j}}{w'_j}$
\end{center}
where each $v_{i,j}\in J\backslash\{0\}$, $u'_{i,j}\in I_1$, and $w'_j\in\mathbb N^+$ does not depend on $i$.

Passing to a subsequence, we may assume that for every $j$, $n_{i,j}$ converges to $\bar n_j\in \bar J=J$, $m_{i,j}$ converges to $\bar m_j\in\bar I_2$, $u_{i,j}$ converges to $\bar u_{j}\in \bar I_1$ and $u'_{i,j}$ converges to $u'_j\in\bar I_1$.

Since $X_i'$ is contained in a bounded family, for every $i$ there exists a very ample Cartier divisor $H_i'$ on $X_i'$ such that $vol(H_i')$ and $-K_{X_i'}\cdot (H_i')^{d-1}$ are bounded from above. In particular, there exists an integer $U>0$ that only depends on $I_1,I_2,J$, such that for any $i$, $-R_0K_{X'_i}\cdot (H'_i)^{d-1}\leq U$. Moreover, for any integral divisor $D_i'$ on $X_i'$, $D_i'\cdot (H_i')^{d-1}$ is an integer.

Let
\begin{center}
$\mathcal{A}_i=R_0K_{X'_i}\cdot (H'_i)^{d-1}+\sum_{j=1}^{p_b}\frac{R_0}{w_j}(w_j-1)(B_{i,j}'\cdot (H'_i)^{d-1})+\sum_{j=1}^{p_c}\frac{R_0}{w_j'}(w_j'-1)(C_{i,j}'\cdot (H'_i)^{d-1})$;

$\mathcal{B}_i=\sum_{j=1}^{p_b}\frac{R_0}{w_j}u_{i,j}(B_{i,j}'\cdot (H'_i)^{d-1})+\sum_{j=1}^{p_c}\frac{R_0}{w'_j}u'_{i,j}(C_{i,j}'\cdot (H'_i)^{d-1})+\sum_{j=1}^{p_m}R_0m_{i,j}(M_{i,j}'\cdot (H'_i)^{d-1})$

$\mathcal{C}_i=\sum_{j=1}^{p_c}\frac{R_0}{w_j'}v_{i,j}(C_{i,j}'\cdot (H'_i)^{d-1})+\sum_{j=1}^{p_n}R_0n_{i,j}(N_{i,j}'\cdot (H_i')^{d-1})$.
\end{center}

Then for any $i$,
\begin{center}
$\mathcal{A}_i+\mathcal{B}_i+c_i\mathcal{C}_i=0$;

$\mathcal{A}_i\geq -U$, $\mathcal{B}_i,\mathcal{C}_i\geq0$.
\end{center}

We state a lemma that will be used in the proof of \hyperref[0105]{Theorem 1.5}.

\begin{lem}\label{0506}\rm
If $1$ is the only accumulation point of $I_1$, and $I_2, J$ has no accumulation points except $\infty$, then \hyperref[Case4]{Case 4} cannot happen.
\end{lem}

\begin{flushleft}\textit{Proof of \hyperref[0506]{Lemma 5.6}}. Since $I_2, J$ has no accumulation points except $\infty$, either $c_i$ converges to $0$, which we have already excluded this case in \textit{Step 1}, or after possibly passing to a subsequence, for every $j$, $n_{i,j}$, $v_{i,j}$, $m_{i,j}$ are all constants. Notice that by our assumption, $u_{i,j}, u'_{i,j}<1-\epsilon$, hence after possibly passing to a subsequence, we may assume that $u_{i,j}, u'_{i,j}$ are all constants. From our construction, $\mathcal{A}_i\in\mathbb Z$, thus $\mathcal{A}_i$ is contained in a finite set that only depends on $U$, and we may assume that $\mathcal{A}_i$ is a constant. But now
\begin{center}
$c_i=\frac{-\mathcal{A_i}-\mathcal{B_i}}{\mathcal{C_i}}$
\end{center}
must be a constant, which is not possible.\end{flushleft}\hfill$\Box$

\begin{flushleft}\textit{Step 7 continued}. Notice that from our construction, $\mathcal{A}_i\in\mathbb Z$, $\mathcal{C}_i\in J$, and $\mathcal{B}_i$ is a sum of elements contained in $I_1,I_2$. Since $\mathcal{A}_i\geq -U$, $\mathcal{B}_i\leq U$, hence $\mathcal{B}_i$ is contained in a DCC set that only depends on $I_1,I_2$ and $U$. 
\end{flushleft}

Since $\mathcal{A}_i\in\mathbb Z$, $\mathcal{A}_i$ is contained in a finite set that only depends on $U$. By our assumption, $\mathcal{C}_i>0$, hence $\mathcal{A}_i<0$. Thus, $\frac{\mathcal{B}_i}{-\mathcal{A}_i}$ is contained in the DCC set 
\begin{center}
$\mathcal{I}:=(\frac{1}{U!}(I_1\cup I_2))_+$,
\end{center}
and $\frac{\mathcal{C}_i}{-\mathcal{A}_i}$ is contained in the DCC set
\begin{center}
$\mathcal{J}=\frac{1}{U!}J$
\end{center}

The accumulation points of $c_i$ are then of the form
\begin{center}
$\frac{1-\alpha}{\beta}$
\end{center}
where $\beta\in\bar{\mathcal{J}}$ and $\alpha\in\bar{\mathcal{I}}$.

Thus $c\in GLCT_{1}(I_1'=\bar{\mathcal{I}},\{0\};J'=\bar{\mathcal{J}})$ by \hyperref[0401]{Lemma 4.1}, hence \hyperref[Case4]{Case 4} holds, and hence \hyperref[0501]{Claim 5.1} holds.\hfill$\Box$

\section{Proof of the theorems}

\textit{Proof of \hyperref[0106]{Theorem 1.6}}. This is immediately followed by \hyperref[0501]{Claim 5.1} and \hyperref[0403]{Lemma 4.3} when $d\geq 3$, \hyperref[0502]{Lemma 5.2} when $d=2$.\hfill$\Box$

\begin{flushleft}\textit{Proof of \hyperref[0105]{Theorem 1.5}}. Let $\tilde J=\{0\}\cup\{\sum_k j_k|j_k\in J\}$. Then $\tilde J\subset [0,\infty)$ is DCC and closed under addition. Notice that \end{flushleft}

(i) If $I,J\subset\mathbb Q$ and all the accumulation points of $I,J$ are contained in $\mathbb Q$, then $(\{1\}\cup I)_+$, $\tilde J$ and all their accumulation points are both contained in $\mathbb Q$.

(ii) If $I$ only has finitely many accumulation points,  $(\{1\}\cup I)_+$ only has finitely many accumulation points.

(iii) If $J$ does not have any accumulation point except $\infty$, $\tilde J$ does not have any accumulation point except $\infty$.

Thus, after possibly replacing $I$ by  $(\{1\}\cup I)_+$, $J$ by $\tilde J$, we may assume $1\in I=I_+$, and $J$ is closed under addition.

Let $I_1=I_2=I$, then the main part of \hyperref[0105]{Theorem 1.5} is followed by \hyperref[0106]{Theorem 1.6} by noticing that $GLCT_{d-1}(I_1',\{0\}; J')\subset GLCT_{d-1}(I_1',I_1'; J')$. By \hyperref[0106]{Theorem 1.6},  from the construction of $J'$ and $I_1'$, if $\bar I_1,\bar I_2$ and $\bar J$ are all contained in $\mathbb Q$, then $GLCT_{d-1}(I_1',I_1'; J')\subset\mathbb Q$, and hence \hyperref[0105]{Theorem 1.5(1)} holds.

To prove (2), by using induction on dimension, we only need to prove that $GLCT_{1}(I_1',\{0\},J')$ only has finitely many accumulation points.

Notice that if $J$ has finitely many accumulation points, since it is closed under addition, $J$ only has $\infty$ as its accumulation point. Moreover, if $I_1,I_2$ has finitely many accumulation points, 
\begin{center}
$\frac{1}{N}(\bar I_1\cup\bar I_2)$
\end{center}
only has finitely many accumulation points. Moreover, let $c>0$ be a lower bound of elements $\not=0$ that contained in $I_1'$. Clearly $c$ only depends on $I_1,I_2$. 

For any $c\in GLCT_{1}(I_1',\{0\},J')$, each $c$ is of the form
\begin{center}
$c=\frac{1-\sum_{l=1}^{p}i_l}{j}$
\end{center}
where $p\leq \frac{1}{c}$ and $i_l\in\frac{1}{N}(\bar I_1\cup\bar I_2)$.

Now we use induction on the order of accumulation points: If all the possible $m$-th order accumulation points of $GLCT_{1}(I_1',\{0\},J')$ are either $0$ or of the form
\begin{center}
$c^m=\frac{1-\sum_{l=1}^{p-m}i_l+\sum_{l=1}^m\lambda_l}{j}$
\end{center}
where $j\in J$, $p\leq\frac{1}{c}$, each $i_l\in\frac{1}{N}(\bar I_1\cup\bar I_2)$ and each $\lambda_l\in\partial(\frac{1}{N}(\bar I_1\cup\bar I_2))$, then the $(m+1)$-th order of $GLCT_{1}(I_1',\{0\},J')$ are either $0$ or of the form
\begin{center}
$c^{m+1}=\frac{1-\sum_{l=1}^{p-m-1}i_l+\sum_{l=1}^{m+1}\lambda_l}{j}$.
\end{center}

In particular, let $m=\lceil\frac{1}{c}\rceil+1$. Then the only $m$-th order accumulation point of $GLCT_{1}(I_1',\{0\},J')$ is $0$, and there's not $(m+1)$-th order accumulation point of  $GLCT_{1}(I_1',\{0\},J')$.\hfill$\Box$

\begin{flushleft}\textit{Proof of \hyperref[0107]{Theorem 1.7}}. 

\textit{Step 1}.
First we show that 
\begin{center}
$\partial GLCT_{d}(I_1,I_2;J)\subset GLCT_{d-1}(I_1,I_2;J)$.
\end{center}\end{flushleft}

By \hyperref[0506]{Lemma 5.6} for dimension $\leq d$, \hyperref[Case4]{Case 4} of \hyperref[0501]{Claim 5.1} cannot happen in any dimension $\leq d$. Now the main part of the theorem follows from \hyperref[0501]{Claim 5.1, Case 1,2,3} for $d\geq 3$ and \hyperref[0504]{Lemma 5.2} for $d=2$.

(1) and (3) are immediately followed from the main part of the theorem. To prove (2), the set of all $(d-2)$-th order accumulation point are contained in $GLCT_{2}(I_1,I_2,J)$, which is equal to $\mathcal{N}_{1}(I_1,I_2;J)$. For any $c\in\mathcal{N}_{1}(I_1,I_2;J)$, $c$ is of the form
\begin{center}
$\frac{2-i_1-i_2}{j}$
\end{center}
where $i_1\in (I_1)_+=I_1$, $i_2\in (I_2)_+=(I_1)_+=I_1$, and $j\in J_+\subset J$. Since $1$ is the only accumulation point of $I_1$, the accumulation points of $\mathcal{N}_{1}(I_1,I_2;J)$ are either $0$, or of the form
\begin{center}
$\frac{1-i_1}{j}$
\end{center}
where $i_1\in I_1$ and $j\in J$. Notice that $\{\frac{1-i_1}{j}|i_1\in I_1, j\in J\}=GLCT_1(I_1,I_2;J)$, and the only possible accumulation point for numbers of the form $\frac{1-i_1}{j}$ is $0$. Thus, the $d$-th order accumulation point set of $GLCT_{d}(I_1,I_2;J)$ is $\{0\}$, and there's no $(d+1)$-th order accumulation point.

\begin{flushleft}\textit{Step 2}. We show that

\begin{center}
$\partial GLCT_{d}(I_1,I_2;J)\supset GLCT_{d-1}(I_1,I_2;J)$,
\end{center}\end{flushleft}
thus finishing the proof.

Suppose $s=glct(X',B'+M',C'+N')$ for some generalized pair $(X',B'+M')$ where $X'$ is of dimension $d-1$, $B'\in I_1$, $M'\in\mathcal{NEF}_{X/X'}(I_2)$, $C'\in J$ and $N'\in\mathcal{NEF}_{X/X'}(J)$ for some $X$ with morphism $g: X\rightarrow X'$. After possibly taking a generalized DLT modification of $(X',B'+M',C'+N')$, we may assume that $X'$ is $\mathbb Q$-factorial. Since $s$ is the generalized log canonical threshold, there exists a generalized log canonical center $F'\subset X'$ of $(X',B'+sC'+M'+sN')$, such that $ga(X',B'+(s+\epsilon)C'+M'+(s+\epsilon)N')<0$ for any $\epsilon>0$.

If $F'$ is not a closed point, after possibly cutting $X'$ by general hyperplanes, we have $s\in GLCT_{d-1-k}(I_1,I_2;J)$ for some $k>0$. Thus, by using induction on dimension and by \hyperref[0401]{Lemma 4.1}, we have
\begin{center}
$s\in GLCT_{d-1-k}(I_1,I_2;J)\subset \partial GLCT_{d-k}(I_1,I_2;J)\subset \partial GLCT_{d}(I_1,I_2;J)$.
\end{center}

Thus, we may assume that $F'=\{p\}$ is a closed point. 

Write $C'=\sum_i c_iC_i'$ into its irreducible components. After possibly taking repeated cyclic covers of $K_X$ and each $C_i'$ near $p$ where we let $f':\bar X'\rightarrow X'$ be the morphism, we may assume that $K_{\bar X'}$ and each $\bar C_i':=f'^*C_i'$ is Cartier near $p$.

Let $\bar B', \bar C',\bar M',\bar N'$ be the pullback of $B',C',M',N'$ on $\bar X'$, $\bar X=\bar X'\times_{X'} X$, $f: \bar X\rightarrow X$ and $\bar g: \bar X\rightarrow \bar X'$ be the corresponding morphisms, and let $\bar M$, $\bar N$ be the pullback of $M,N$ on $\bar X$. For any prime divisor $E$ on $X$, let $\bar E$ be the prime divisor on $\bar X$ that dominates $E$. 

Near the generic point of $\bar E$, for any generalized pair $(X',\Delta'+\Gamma')$ with data $(\Gamma, X\rightarrow X')$, we let $\bar\Delta'=f'^*\Delta'$ and $\bar\Gamma'=f'^*\Gamma'$. Suppose $r$ is the ramification index of $f$ along $\bar E$.

Then since $f'$ is unramified in codimension $1$, we have

\begin{center}
$K_{\bar X}=\bar g^*(K_{\bar X'})+(ga(\bar E,\bar X',\bar\Delta'+\bar\Gamma')-1)\bar E=\bar g^*f'^*K_{X'}+(ga(\bar E,\bar X',\bar\Delta'+\bar\Gamma')-1)\bar E$
\end{center}
and

\begin{center}
$K_{\bar X}=f^*K_{X}+(r-1)\bar E=f^*(g^*(K_{X'}+(ga(E,X',\Delta'+\Gamma')-1)E)+(r-1)\bar E=f^*g^*K_{X'}+(r\cdot ga(E,X',\Delta'+\Gamma')-1)\bar E$.
\end{center}

Thus, since $\bar g^*f'^*K_{X'}=f^*g^*K_{X'}$, we have $r\cdot ga(E,X',\Delta'+\Gamma')=ga(\bar E,\bar X',\bar\Delta'+\bar\Gamma')$. 

In particular, we have $s=glct(\bar X',\bar B'+\bar M';\bar C'+\bar N')$. Thus we may replace $X',X,B',C',M',N',M,N$ by $\bar X',\bar X,\bar B',$ $\bar C',\bar M',\bar N',\bar M,\bar N$, and assume that each component of $C_i'$ is Cartier near $p$. 

After possibly replacing $X'$ by a neighborhood near $p$ and possibly replacing $X$, we may assume that $X'\cong Spec\ R$. Now consider the Cartier divisors $D_{i,k}'$ contained in $Spec (R,x_d)\cong X'\times\mathbb A^1$ that corresponds to the ideal $(\mathfrak{c}_i,x_{d}^k)$ for every integer $k>0$, where $\mathfrak{c}_i$ is defined by
\begin{center}
$\mathfrak{c}_i\cdot\mathcal{O}_{X}=\mathcal{O}_X(-C_i')$.
\end{center}

We consider 
\begin{center}
$s_{k}:=glct(X'\times\mathbb A^1, B'\times\mathbb A^1+M'\times\mathbb A^1,\sum_i c_iD_{i,k}'+N_i'\times\mathbb A^1)$
\end{center}

Then since each $s_i$ is a fixed positive real number, by following the same lines of the proof of \hyperref[Kol97]{[Kol97, Proposition 8.21]}, for $k\gg 0$, we have $s_k=s+\frac{1}{k}$, thus after possibly passing to a subsequence, $s_k$ forms a strict decreasing sequence that converges to $s$. Since $s_k\in GLCT_{d}(I_1,I_2;J)$ by our construction, we conclude that $\partial GLCT_{d}(I_1,I_2;J)\supset GLCT_{d-1}(I_1,I_2;J)$, and the proof is finished.

\hfill$\Box$
\section{References}

\hspace{6mm}[Bir07]\label{Bir07} Caucher Birkar, \textit{Ascending chain condition for log canonical thresholds and termination of log flips.} Duke Math. J. Volume \textbf{136}, Number 1 (2007), 173-180. 

[Bir16]\label{Bir16} Caucher Birkar, \textit{Singularities of linear systems and boundedness of Fano varieties.} arXiv: 1609.05543v1.\vspace{2mm}

[BS05]\label{BS05} Caucher Birkar and Vyacheslav Vladimirovich Shokurov, \textit{Mld's vs thresholds and flips}. J. Reine Angew. Math. \textbf{638} (2010), 209-234.\vspace{2mm}

[BZ16]\label{BZ16} Caucher Birkar and De-Qi Zhang, \textit{Effectivity of Iitaka fibrations and pluricanonical systems of polarized pairs}. Publ. Math. Inst. Hautes Etudes Sci., 123:283-331, 2016.\vspace{2mm}

[HMX14]\label{HMX14} Christopher D.Hacon, James M\textsuperscript{c}Kernan, and Chenyang Xu, \textit{ACC for log canonical thresholds}, Ann. of Math. (2) \textbf{180} (2014), no. 2, pp.523-571.\vspace{2mm}

[KM98]\label{KM98} J\'{a}nos Koll\'{a}r and Shigefumi Mori, \textit{Birational geometry of algebraic varieties}, Cambridge Tracts in Math. \textbf{134}, Cambridge Univ. Press, (1998).\vspace{2mm}

[Kol97]\label{Kol97} J\'{a}nos Koll\'{a}r, \textit{Singularities of pairs}, Algebraic  geometry (Santa Cruz,1995), 221–287. Proc. Symp. Pure Math. 62, Amer. Math. Soc. (1997).\vspace{2mm}

[Kol08]\label{Kol08} J\'{a}nos Koll\'{a}r, \textit{Which powers of holomorphic functions are integrable?}. arXiv:0805.0756.\vspace{2mm}

[Sho93]\label{Sho93} Vyacheslav Vladimirovich Shokurov, \textit{3-fold log flips}, With an appendix by Yujiro Kawamata, Russian Acad. Sci. Izv. Math. \textbf{40} (1993), no. 1, 95-202.\vspace{2mm}

[Tot10]\label{Tot10} Burt Totaro, \textit{The ACC Conjecture for log canonical thresholds [after de Fernex, Ein, Musta\c{t}\u{a}, Koll\'{a}r]}, S\'{e}minaire BOURBAKI, 62\`{e}me ann\'{e}e, 2009-2010, no 1025, Juin 2010, 15p.
\end{document}